\documentclass[12pt]{article}
\usepackage{amsmath}

\usepackage{amssymb}
\usepackage{graphicx}
\usepackage{epsfig}

\begin{document}

{\large Identification of a diffusion coefficient in strongly degenerate
parabolic equations with interior degeneracy}

\bigskip

Genni Fragnelli$^{1},$ Gabriela Marinoschi$^{2},$ Rosa Maria Mininni$^{1}$
and Silvia Romanelli$^{1}$

\bigskip

$^{1}${\small Dipartimento di Matematica, Universit\`{a} degli Studi di Bari
Aldo Moro, via E. Orabona 4, 70125 Bari (Italy)}

$^{2}${\small Institute of Mathematical Statistics and Applied Mathematics
of the Romanian Academy, Calea 13 Septembrie 13, Bucharest (Romania)}

{\small \bigskip }

{\small E-mail: genni.fragnelli@uniba.it, gmarino@acad.ro,
rosamaria.mininni@uniba.it and silvia.romanelli@uniba.it}

\bigskip

\textbf{Abstract}

We study two identification problems in relation with a strongly degenerate
parabolic diffusion equation characterized by a vanishing diffusion
coefficient $u\in W^{1,\infty },$ with the property $\frac{1}{u}\notin
L^{1}. $ The aim is to identify $u$ from certain observations on the
solution, by a technique of nonlinear optimal control with control in
coefficients. The existence of a controller $u$ which is searched in $%
W^{1,\infty }$ and the determination of the optimality conditions are given
for homogeneous Dirichlet boundary conditions. An approximating problem
further introduced allows a better characterization of the optimality
conditions, due to the supplementary regularity of the approximating state
and dual functions and to a convergence result. Finally, an identification
problem with final time observation and homogeneous Dirichlet-Neumann
boundary conditions in the state system is considered. By using more
technical arguments we provide the explicit form of $u$ and its uniqueness.

\bigskip

\noindent Keywords: inverse problems, degenerate diffusion equations, internal degeneracy, optimal control, optimality conditions

\bigskip

\noindent MSC 2010: 35K65, 35R30, 49N45, 49KXX

\bigskip

\section{Introduction}

In this article we study two identification problems in relation with an
evolution equation with a second order differential degenerate operator, in
divergence form, $A_{1}y:=(uy_{x})_{x}$, when the diffusion coefficient $u$
vanishes at an interior point of an one dimensional space domain. The
purpose is to determine $u$ from certain observations on the solution to the
evolution equation.

Degenerate parabolic operators naturally arise in many problems:
Budyko-Sellers models in climatology (see, e.g., \cite{tv}), boundary layer
models in physics (see, e.g., \cite{br}), Wright-Fisher and Fleming-Viot
models in genetics (see, e.g., \cite{fleming-viot-79}, \cite{s}),
Black-Merton-Scholes models in mathematical finance (see, e.g., \cite{mro}).
The degenerate operator $A_{1}$ has been studied under different boundary
conditions, see, for example, \cite{cmp}, \cite{cty}, \cite{gl}, \cite{v}.
In \cite{fav-yagi}, \cite{gl}, the authors consider degenerate operators
with boundary conditions of Dirichlet, Neumann, periodic, or nonlinear Robin
type. In \cite{acf} the authors consider the degenerate operator in
divergence and in non divergence form with Dirichlet or Neumann boundary
conditions, giving more importance to controllability problems of the
associated parabolic evolution equations. However, all previous papers deal
with a degenerate operator with degeneracy at the boundary of the domain,
for example of the form of the double power function $%
u(x)=x^{k_{1}}(1-x)^{k_{2}},$ $x\,\in \lbrack 0,1],$ where $k_{1}$ and $%
k_{2} $ are positive constants.

To the best of our knowledge, Stahel's paper \cite{s1} is the first treating
a problem with a degeneracy which may be interior. In particular, Stahel
considers a parabolic problem in ${\mathbb{R}}^{N}$ with Dirichlet, Neumann
or mixed boundary conditions, associated with a $N\times N$ matrix $a$
(defining the diffusion coefficient), which is positive definite and
symmetric, but whose smallest eigenvalue might converge to $0$ as the space
variable approaches a singular set contained in the closure of the space
domain. In this case, he proves that the corresponding abstract Cauchy
problem has a solution, provided that $\underline{a}^{-1}\in L^{q}(\Omega ,{%
\mathbb{R}})$ for some $q>1$, where $\Omega $ is an open bounded subset of $%
\mathbb{R}^{N},$ and $\underline{a}(x):=\min \{a(x)\xi \cdot \xi :\Vert \xi
\Vert =1\}.$ In \cite{cm} there is treated a class of variational degenerate
elliptic problems with interior and boundary degeneracy in the case that
there exists $k\in (0,2]$ such that $u$ decreases more slowly than $%
\left\vert x-z\right\vert ^{k}$ near every point $z\in u^{-1}\{0\}$ in the
case of bounded domains. The assumption regarding the interior degeneracy of
$u$ is generalized in \cite{frggr1} when $N=1$. In particular, in \cite%
{frggr1} the authors analyze in detail degenerate operators in divergence
and non-divergence form, under Dirichlet boundary conditions in spaces $%
L^{2}(0,1)$ with or without weight and show that under suitable assumptions
they generate analytic semigroups. In \cite{fm} the same operators with an
interior degeneracy are considered, giving more importance to the null
controllability of the associated parabolic equations. In \cite%
{af-gm-identif} a control problem involving a nonlinear nonautonomous
operator degenerating on a positive measure interior subset of the space
domain is studied.

The interest in this kind of interior degeneracy problems is due to the fact
that they govern diffusion of a substance in water, soil or air, heat flow
in a material, diffusion of a population in a habitat. The nonhomogeneity of
the medium is expressed by the space dependence of the diffusion coefficient
with its possible vanishing at some points. For example, a certain composite
material can block the heat flow at a certain point, or the migration of
small mammal species can degenerate due to environmental heterogeneity and
barriers (see, for example, \cite{b}, \cite{ly} and the references therein).
For this reason it is important to study identification and control problems
associated to these degenerate equations. In particular, in \cite{b} and
\cite{ly} the authors consider two optimal problems: the first one is to
minimize the damage and trapping costs, the second one is to maximize the
difference between harvesting cost and economic revenue. Another application
can be related to the study of the design of biological channels (see, e.g.,
\cite{gm-metab}) in the case when the metabolite diffusion coefficient
vanishes.

From the mathematical point of view and in connection with the work of
Fragnelli et al. (see \cite{frggr1}) we focus on identifying, on the basis
of some observations, the diffusion coefficient in the degenerate parabolic
equation%
\begin{equation}
\frac{\partial y}{\partial t}-\frac{\partial }{\partial x}\left( u(x)\frac{%
\partial y}{\partial x}\right) =f\mbox{ in }Q:=(0,T)\times (0,L),\mbox{ }%
T,L\in (0,\infty ),  \label{1}
\end{equation}%
with the initial condition
\begin{equation}
y(0,x)=y_{0}(x)\mbox{ in }(0,L),  \label{2}
\end{equation}%
and various boundary conditions, in particular of homogeneous Dirichlet type%
\begin{equation}
y(t,0)=y(t,L)=0\mbox{ in }(0,T).  \label{3}
\end{equation}%
We envisage a vanishing $u$ at a point $x_{0}\in (0,L)$ considered in \cite%
{frggr1} by taking also into account the behavior of the function $\frac{1}{u%
},$ corresponding either to a so-called slow diffusion (the strongly
degenerate case) or to a fast diffusion (the weakly degenerate case).

As far as we know, the particular properties related to the behavior of $%
\frac{1}{u}$ have not been considered in other identification or control
papers.

To motivate our work we specify that the effects of the vanishing diffusion
coefficient upon the solution to the diffusion equation (concentration,
temperature, density) change according to the particular form of $u(x).$
This behavior is illustrated in the graphics realized with Comsol
Multiphysics v. 3.5a (FLN License 1025226), e.g., for $u(x)=\left\vert
x-x_{0}\right\vert ^{n},$ $x\in \lbrack 0,1],$ $x_{0}=0.5,$ $y_{0}=1.$ Fig.
1 represents the values of the solution $y(t,x)$ to (\ref{1})-(\ref{3})
along $Ox,$ for $t=0;$ $0.5;$ $1;$ 1.5; 2$,$ for $n=2,$ and Fig. 2 shows the
distribution of the values of the solution $y(t,x)$ along $Ox$ at the same
times for $n=4.$ We observe that if $u$ has a zero of a higher order of
multiplicity ($n=4)$ at $x_{0}$, then the solution lies at high values in a
larger subset of $(0,1).$

\begin{center}
    \begin{tabular}{cc} 
           \includegraphics[width=2.6in]{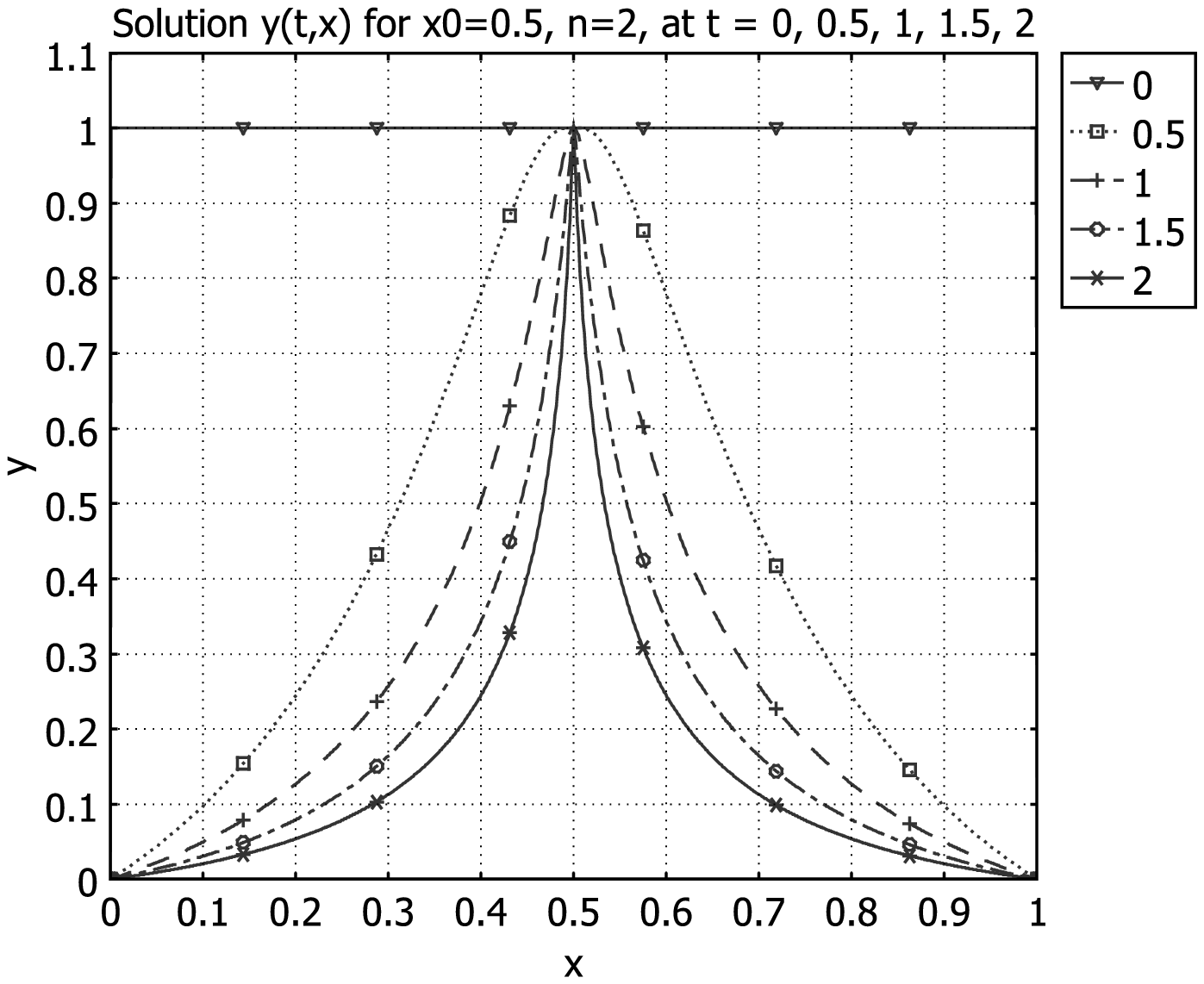} & \includegraphics[width=2.6in]{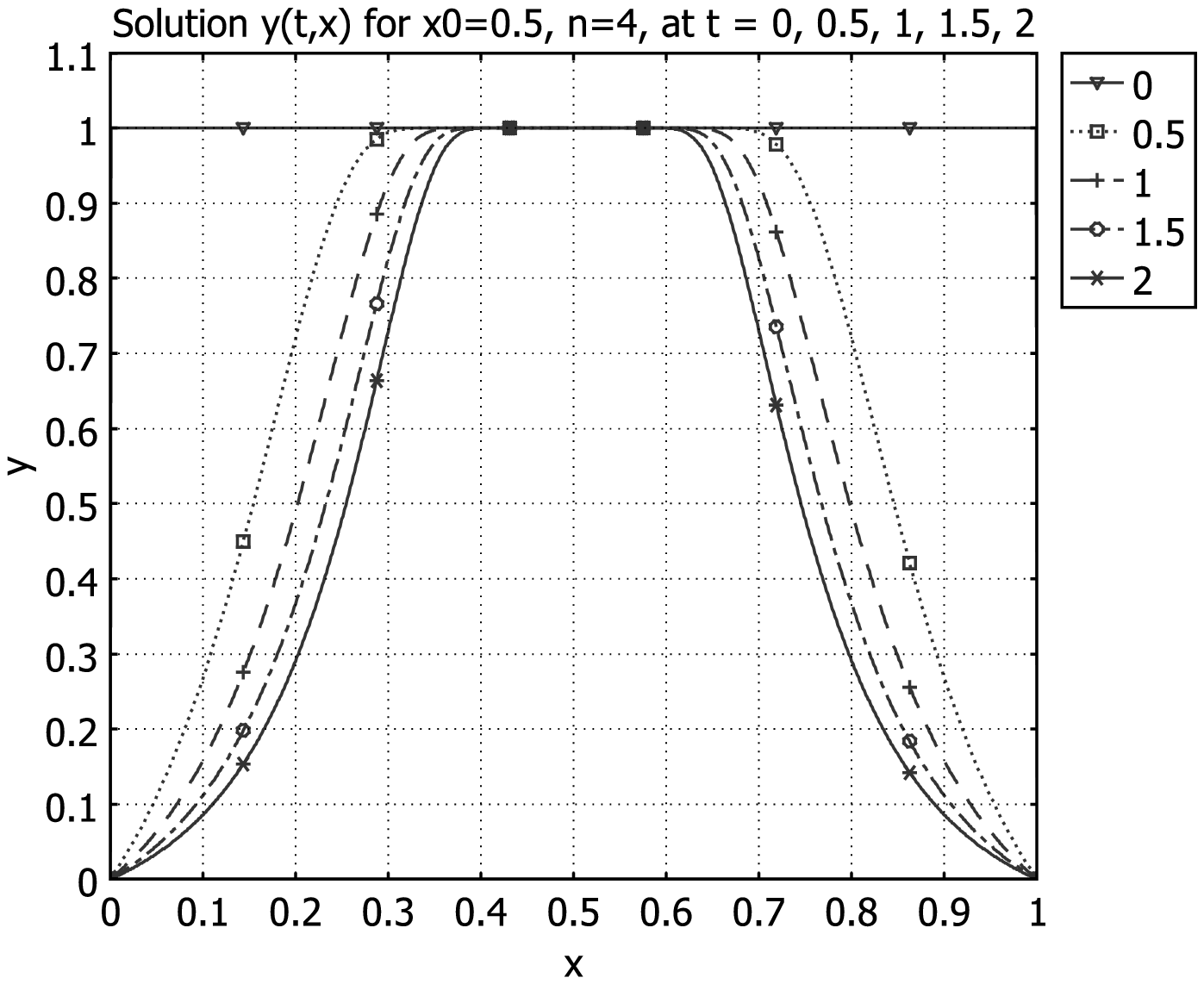} \\
           Fig. 1. $y(t,x)$ for $u(x)=\left\vert x-x_{0}\right\vert ^{2}$ &      Fig. 2. $y(t,x)$ for $u(x)=\left\vert x-x_{0}\right\vert ^{4}$
           \\ 
    \end{tabular}
\end{center}

In some practical situations, for instance in a pollutant diffusion process
or in a diffusion of a population (bacteria e.g.) in a medium, it is of
interest to identify the diffusion coefficient which is suspected to
determine high levels of concentration (or density) $y$ in a certain subset
of the flow domain due to the possible diffusion stopping at some point $%
x_{0}.$ A similar interest is in the design of a composite material for
determining the material properties which preserve the temperature at
certain high values.

We mainly aim to determine the function $u(x)$ in the system (\ref{1})-(\ref%
{3}) from the observation of the spatial mean ($M_{T})$ of the state at a
final time $T.$ However, since other physical quantities, as the mean value $%
M$ or the mean flux $M_{f}$ over $Q$ may be measured, we shall identify $u$
by combining these possible observations in a unique mathematical problem of
minimization of a cost functional. Namely, first we study the problem%
\begin{align*}
\mbox{Minimize }& \left\{ \frac{\lambda _{1}}{2}\left(
\int_{Q}u(x)y_{x}^{u}(t,x)dxdt-M_{f}\right) ^{2}+\frac{\lambda _{2}}{2}%
\left( \int_{0}^{L}y^{u}(T,x)dx-M_{T}\right) ^{2}\right. \\
& \left. +\frac{\lambda _{3}}{2}\left( \int_{Q}y^{u}(t,x)dxdt-M\right)
^{2}\right\}
\end{align*}%
subject to (\ref{1})-(\ref{3}), for all $u$ in a set $U,$ which will be
specified in Section 2.2. Here, $M_{f},$ $M_{T},$ $M,$ $\lambda _{1},$ $%
\lambda _{2},$ $\lambda _{3},$ are nonnegative real numbers, and there
exists at least one $i\in \{1,2,3\}$ such that $\lambda _{i}>0.$

The notation $y^{u}$ indicates the solution to (\ref{1})-(\ref{3})
corresponding to $u.$

The various choices of the constants $\lambda _{i},$ $i=1,2,3,$ enhance a
higher or lower importance to the terms in the functional, according to the
problem requirements.

We provide a detailed study of the diffusion equation in the divergence form
(\ref{1}) for the strongly degenerate case (see Definition 2.1) both for
homogeneous Dirichlet boundary conditions and for homogeneous
Dirichlet-Neumann boundary conditions.

We prove the existence of a solution to the above nonlinear minimization
problem, by approaching it as an optimal control problem in coefficients.
The diffusion coefficient $u$ is determined by the necessary conditions that
$u,$ as a solution to $(P),$ should satisfy. They are called optimality
conditions.\ To this end we study the existence in the state system, the
system in variations and the dual system associated to the minimization
problem. Then we present a detailed computation of the necessary conditions
of optimality which is quite technical due to the fact that we search for a
function $u$ in $W^{1,\infty }(0,L).$ In the second identification problem
we require to identify $u$ from the final observation, imposing to
\begin{equation*}
\mbox{Minimize \ }\left\{ \int_{0}^{L}y^{u}(T,x)dx\right\}
\end{equation*}%
for all $u\in U,$ subject to (\ref{1})-(\ref{2}), with homogeneous
Dirichlet-Neumann boundary conditions. In this case we exploit some further
properties of the derivatives of the state in order to get a more explicit
form of $u$ and its uniqueness.

\section{Homogeneous Dirichlet boundary conditions}

\subsection{Preliminaries and the state system}

We begin with some notation, definitions and results given in \cite{frggr1},
considering the operator $A_{1}y=(uy_{x})_{x}$ in the strongly degenerate
case. For simplicity we denote $\phi _{t},$ $\phi _{x}$ the partial
derivatives of a generic function $\phi (t,x)$ with respect to $t$ and $x.$

\medskip

\noindent \textbf{Definition 2.1.} $A_{1}y=(uy_{x})_{x}$ is called \textit{%
strongly degenerate} if there exists $x_{0}\in (0,L)$ such that $u(x_{0})=0,$
$u(x)>0$ on $[0,L]\backslash \{x_{0}\},$ $u\in W^{1,\infty }(0,L)$ and $%
\frac{1}{u}\notin L^{1}(0,1).$

\medskip

As an example we can mention $u(x)=\left\vert x-x_{0}\right\vert ^{k},$ $%
k\geq 1.$

\medskip

We define the weighted space
\begin{align}
H_{u}^{1}(0,L)& =\{y\in L^{2}(0,L);\mbox{ }y\mbox{ locally absolutely
continuous in }[0,L]\backslash \{x_{0}\},  \notag \\
\sqrt{u}y_{x}& \in L^{2}(0,L),\mbox{ }y(0)=y(L)=0\},  \label{4}
\end{align}%
with the norm
\begin{equation}
\left\Vert y\right\Vert _{H_{u}^{1}(0,L)}=\left( \left\Vert y\right\Vert
_{L^{2}(0,L)}^{2}+\left\Vert \sqrt{u}y_{x}\right\Vert
_{L^{2}(0,L)}^{2}\right) ^{1/2}.  \label{5}
\end{equation}%
According to \cite{frggr1}, Proposition 2.3, we have that%
\begin{multline*}
H_{u}^{1}(0,L)=\{y\in L^{2}(0,L);\mbox{ }y\mbox{ locally absolutely
continuous in }[0,L]\backslash \{x_{0}\}, \\
\sqrt{u}y_{x}\in L^{2}(0,L),\mbox{ }uy\mbox{ is continuous at }x_{0},\mbox{ }%
y(0)=y(L)=(uy)(x_{0})=0\}.
\end{multline*}%
We specify that $H_{u}^{1}(0,L)$ is a Hilbert space and $H_{u}^{1}(0,L)%
\hookrightarrow L^{2}(0,L)\hookrightarrow (H_{u}^{1}(0,L))^{\prime },$ where
$(H_{u}^{1}(0,L))^{\prime }$ is the dual of $H_{u}^{1}(0,L)$ and "$%
\hookrightarrow $" means a continuous and dense embedding$.$

For simplicity, we denote
\begin{equation*}
H=L^{2}(0,L),\mbox{ }V_{u}=H_{u}^{1}(0,L),\mbox{ }V_{u}^{\prime
}=(H_{u}^{1}(0,L))^{\prime },
\end{equation*}%
where they indicate subscripts.

Let us consider
\begin{equation}
H_{u}^{2}(0,L)=\{y\in H_{u}^{1}(0,L);\mbox{ }uy_{x}\in H^{1}(0,L)\}
\label{6}
\end{equation}%
and define
\begin{equation*}
A_{1}y:=(uy_{x})_{x},\mbox{ }A_{1}:D(A_{1})\subset L^{2}(0,L)\rightarrow
L^{2}(0,L),\mbox{ }D(A_{1})=H_{u}^{2}(0,L).
\end{equation*}%
According to \cite{frggr1}, Proposition 2.4,
\begin{eqnarray}
D(A_{1}) &=&\{y\in L^{2}(0,L);\mbox{ }y\mbox{ is absolutely continuous in }%
[0,L]\backslash \{x_{0}\},  \label{6-0} \\
uy &\in &H_{0}^{1}(0,L),\mbox{ }uy_{x}\in H^{1}(0,L),\mbox{ }%
uy(x_{0})=uy_{x}(x_{0})=0\}.  \notag
\end{eqnarray}

By \cite{frggr1}, Theorem 2.7, $A_{1}:D(A_{1})\rightarrow L^{2}(0,L)$ is
self-adjoint, nonpositive on $L^{2}(0,L)$ and it generates a positivity
preserving semigroup. This result is used further to prove that (\ref{1})-(%
\ref{3}) has a unique mild solution if $y_{0}\in L^{2}(0,L)$ and $f\in
L^{2}(Q).$ This is a strong solution if $y_{0}\in H_{u}^{2}(0,L)$ and also
if $y_{0}\in H_{u}^{1}(0,L)$ (see Theorem 4.1 and Remark 4.2 in \cite{frggr1}%
).

We also recall the following result (Lemma 2.6 in \cite{frggr1}):

\medskip

\noindent \textbf{Lemma 2.2. }\textit{For all }$(y,z)\in
H_{u}^{2}(0,L)\times H_{u}^{1}(0,L)$\textit{\ one has}%
\begin{equation}
\int_{0}^{L}(uy_{x})_{x}zdx=-\int_{0}^{L}uy_{x}z_{x}dx.  \label{7}
\end{equation}

\medskip

For the purposes of our paper we present the existence result for the state
system (\ref{1})-(\ref{3}) by a variational way. For convenience, and where
no confusion can be made we shall not write the function arguments in the
integrands.

\medskip

\noindent \textbf{Definition 2.3. }Let $y_{0}\in L^{2}(0,L),$ $f\in
L^{2}(0,T;(H_{u}^{1}(0,L))^{\prime })$ and $u$ with the properties of
Definition 2.1. We call a \textit{solution} to (\ref{1})-(\ref{3}) a
function
\begin{equation}
y\in C([0,T];L^{2}(0,L))\cap L^{2}(0,T;H_{u}^{1}(0,L))\cap
W^{1,2}([0,T];(H_{u}^{1}(0,L))^{\prime }),  \label{18-10}
\end{equation}%
which satisfies the equation
\begin{equation}
\int_{0}^{T}\left\langle \frac{dy}{dt}(t),\psi (t)\right\rangle
_{V_{u}^{\prime },V_{u}}dt+\int_{Q}uy_{x}\psi
_{x}dxdt=\int_{0}^{T}\left\langle f(t),\psi (t)\right\rangle _{V_{u}^{\prime
},V_{u}}dt,  \label{18-11}
\end{equation}%
for any $\psi \in L^{2}(0,T;H_{u}^{1}(0,L)),$ and the initial condition $%
y(0)=y_{0}.$

\medskip

\noindent \textbf{Theorem 2.4.} \textit{If }$y_{0}\in L^{2}(0,L),$ $f\in
L^{2}(0,T;(H_{u}^{1}(0,L))^{\prime }),$ \textit{then }(\ref{1})-(\ref{3})
\textit{has a unique solution }%
\begin{equation*}
y\in C([0,T];L^{2}(0,L))\cap L^{2}(0,T;H_{u}^{1}(0,L))\cap
W^{1,2}([0,T];(H_{u}^{1}(0,L))^{\prime }),
\end{equation*}%
\textit{satisfying the estimate}%
\begin{equation}
\sup_{t\in \lbrack 0,T]}\left\Vert y(t)\right\Vert
_{H}^{2}+\int_{0}^{T}\left\Vert y(t)\right\Vert _{V_{u}}^{2}dt\leq
C_{T}(\left\Vert y_{0}\right\Vert _{H}^{2}+\left\Vert f\right\Vert
_{L^{2}(0,T;V_{u}^{\prime })}^{2}).  \label{7-1}
\end{equation}%
\textit{If, in addition }$y_{0}\in H_{u}^{1}(0,L)$\textit{\ and }$f\in
L^{2}(Q),$ \textit{then}%
\begin{equation}
y\in W^{1,2}([0,T];L^{2}(0,L))\cap L^{2}(0,T;H_{u}^{2}(0,L))\cap L^{\infty
}(0,T;H_{u}^{1}(0,L))  \label{7-1-0}
\end{equation}%
\textit{and it satisfies}
\begin{equation}
\sup_{t\in \lbrack 0,T]}\left\Vert y(t)\right\Vert
_{V_{u}}^{2}+\int_{0}^{T}\left( \left\Vert \frac{dy}{dt}(t)\right\Vert
_{H}^{2}+\left\Vert (uy_{x})_{x}(t)\right\Vert _{H}^{2}\right) dt\leq
C_{T}\left( \left\Vert y_{0}\right\Vert _{V_{u}}^{2}+\left\Vert f\right\Vert
_{L^{2}(Q)}^{2}\right) ,  \label{7-2}
\end{equation}%
\textit{where }$C_{T}$\textit{\ denotes several positive constants.}

\textit{If} $y_{0}\geq 0$\textit{\ a.e. on }$(0,L)$ \textit{and }$f\geq 0$
\textit{a.e. on} $Q,$ \textit{then }$y(t)\geq 0$\textit{\ a.e. in }$(0,L)$%
\textit{\ for all }$t\in \lbrack 0,T].$

\medskip

\noindent \textbf{Proof. }Let us introduce the linear operator
\begin{equation*}
A:H_{u}^{1}(0,L)\rightarrow (H_{u}^{1}(0,L))^{\prime }
\end{equation*}%
by
\begin{equation}
\left\langle Az,\psi \right\rangle _{V_{u}^{\prime
},V_{u}}=\int_{0}^{L}u(x)z_{x}(x)\psi _{x}(x)dx,\mbox{ for any }\psi \in
H_{u}^{1}(0,L).  \label{18-00}
\end{equation}%
It is continuous and monotone%
\begin{equation}
\left\Vert Az\right\Vert _{V_{u}^{\prime }}=\sup_{\psi \in V_{u},\left\Vert
\psi \right\Vert _{V_{u}}\leq 1}\left\vert \left\langle Az,\psi
\right\rangle _{V_{u}^{\prime },V_{u}}\right\vert \leq \left\Vert
z\right\Vert _{V_{u}},  \label{18-01}
\end{equation}%
\begin{equation}
\left\langle Az,z\right\rangle _{V_{u}^{\prime },V_{u}}\geq 0,
\label{18-01-0}
\end{equation}%
and has the property%
\begin{equation}
\left\langle Az,z\right\rangle _{V_{u}^{\prime
},V_{u}}=\int_{0}^{L}uz_{x}^{2}dx=\left\Vert z\right\Vert
_{V_{u}}^{2}-\left\Vert z\right\Vert _{H}^{2}.  \label{18-02}
\end{equation}%
We consider the Cauchy problem
\begin{eqnarray}
\frac{dy}{dt}(t)+Ay(t) &=&f(t),\mbox{ a.e. }t\in (0,T)  \label{18-03} \\
y(0) &=&y_{0}.  \notag
\end{eqnarray}%
If%
\begin{equation*}
y_{0}\in L^{2}(0,L),\mbox{ }f\in L^{2}(0,T;(H_{u}^{1}(0,L))^{\prime })
\end{equation*}%
then the Cauchy problem has a unique solution belonging to the spaces
indicated in (\ref{18-10}) (see e.g. \cite{lions}, p. 162). Obviously, this
is a solution to (\ref{1})-(\ref{3}) in the sense of Definition 2.3, so that
(\ref{18-11}) and (\ref{18-03}) are equivalent. The estimate (\ref{7-1})
follows by setting in (\ref{18-11}) $\psi =y$ \ and performing a few
computations involving the Gronwall's lemma.

Now, we observe that $-A_{1}$ is the restriction of $A$ on $L^{2}(0,L).$
Indeed, if $y\in H_{u}^{2}(0,L),$ by (\ref{7}) we see that $-A_{1}$
coincides with $A.$ Therefore, if $y_{0}\in D(A_{1})=H_{u}^{2}(0,L)$ and $%
f\in W^{1,1}([0,T];L^{2}(0,L)),$ problem (\ref{18-03}) with $A$ replaced by
its restriction $-A_{1}$ has a more regular solution $y\in W^{1,\infty
}([0,T];L^{2}(0,L))\cap L^{\infty }(0,T;H_{u}^{2}(0,L))$. This follows by
Theorem 4.9 (see \cite{vb-springer-2010}, p. 151). If $y_{0}\in
H_{u}^{1}(0,L)$ and $f\in L^{2}(Q),$ the density of $H_{u}^{2}(0,L)$ in $%
H_{u}^{1}(0,L),$ the density of $W^{1,1}([0,T];L^{2}(0,L))$ in $L^{2}(Q)$
and some standard estimates lead to the second part of the theorem.

Finally, the nonnegativity of the solution follows by the positivity
preserving property of the semigroup generated by $A_{1}.$\hfill $\square $

\subsection{Existence in $(P)$}

We denote
\begin{eqnarray}
J(u) &=&\frac{\lambda _{1}}{2}\left( \int_{Q}uy_{x}^{u}dxdt-M_{f}\right)
^{2}+\frac{\lambda _{2}}{2}\left( \int_{0}^{L}y^{u}(T,x)dx-M_{T}\right) ^{2}
\label{J} \\
&&+\frac{\lambda _{3}}{2}\left( \int_{Q}y^{u}dxdt-M\right) ^{2}  \notag
\end{eqnarray}%
and introduce the minimization problem
\begin{equation}
\mbox{Minimize }J(u)\mbox{ for all }u\in U,  \tag{$P$}
\end{equation}%
subject to (\ref{1})-(\ref{3}), where%
\begin{eqnarray}
U &=&\{u\in W^{1,\infty }(0,L);\mbox{ }u_{m}(x)\leq u(x)\leq u_{M}(x),\mbox{
}  \label{U} \\
u(0) &=&u_{0},\mbox{ }u(L)=u_{L},\mbox{ }\left\vert u_{x}(x)\right\vert \leq
u_{\infty }\mbox{ a.e. }x\in (0,L)\}.  \notag
\end{eqnarray}%
We assume the following hypotheses:
\begin{eqnarray}
x_{0} &\in &(0,L),\mbox{ }u_{\infty }\in \lbrack 0,\infty ),  \label{3-1} \\
u_{m},\mbox{ }u_{M} &\in &C[0,L],\mbox{ }u_{M}(x)\leq \alpha (x)u_{m}(x)%
\mbox{ for }x\in \lbrack 0,L],\mbox{ }\alpha \in C[0,L],\alpha \geq 1,
\notag \\
0 &<&u_{m}(x)<u_{M}(x)\mbox{ for }x\in \lbrack 0,L]\backslash \{x_{0}\},%
\mbox{ }u_{m}(x_{0})=u_{M}(x_{0})=0,  \notag
\end{eqnarray}%
and
\begin{equation}
\int_{0}^{L}\frac{1}{u_{M}(x)}dx=+\infty .  \label{3-2}
\end{equation}%
Then, by (\ref{3-1}), for all $u\in U$ we have
\begin{equation}
0<u_{m}(0)\leq u_{0}\leq u_{M}(0),\mbox{ }0<u_{m}(L)\leq u_{L}\leq u_{M}(L).
\label{3-2-0}
\end{equation}%
These conditions ensure the fact that the operator $A_{1}$ is strongly
degenerate, because $u\in U$ and (\ref{3-1}) imply that $u(x_{0})=0$, $u>0$
in $[0,L]\backslash \{x_{0}\},$ and (\ref{3-2}) establishes that $\frac{1}{u}%
\notin L^{1}(0,L).$

Moreover, the assumption $u_{M}(x)\leq \alpha (x)u_{m}(x)$ implies that if $%
u,v\in U$ then $\left\vert \frac{v}{u}(x)\right\vert \leq \left\Vert \alpha
\right\Vert _{L^{\infty }(0,L)}$ for $x\in \lbrack 0,L]\backslash \{x_{0}\}$
and
\begin{equation}
H_{u}^{1}(0,L)=H_{v}^{1}(0,L)\mbox{ for any }u,v\in U.  \label{3-3}
\end{equation}%
Indeed, if $y\in H_{u}^{1}(0,L)$ then $y$ is locally absolutely continuous
in $[0,L]\backslash \{x_{0}\},$ $y\in L^{2}(0,L),$ $\sqrt{u}y_{x}\in
L^{2}(0,L),$ $y(0)=y(L)=0.$ Let $v\in U$. Then, by a simple calculation
\begin{equation}
\int_{0}^{L}vy_{x}^{2}dx\leq \left\Vert \alpha \right\Vert _{L^{\infty
}(0,L)}\left\Vert \sqrt{u}y_{x}\right\Vert _{L^{2}(0,L)}^{2}  \label{3-3-0}
\end{equation}%
and so $y\in H_{v}^{1}(0,L).$ Analogously, we get $H_{v}^{1}(0,L)\subset
H_{u}^{1}(0,L).$

\medskip

\noindent \textbf{Theorem 2.5. }\textit{Let }$y_{0}\in L^{2}(0,L),$\textit{\
}$y_{0}\geq 0$ \textit{on }$(0,L)$\textit{, }$f\in L^{2}(Q),$ $f\geq 0$%
\textit{\ a.e. on} $Q.$\textit{\ Then }$(P)$\textit{\ has at least one
solution }$u$ \textit{with the corresponding state}
\begin{equation*}
y\in C([0,T];L^{2}(0,L))\cap L^{2}(0,T;H_{u}^{1}(0,L))\cap
W^{1,2}([0,T];(H_{u}^{1}(0,L))^{\prime }).
\end{equation*}%
\textit{If }$y_{0}\in H_{u}^{1}(0,L),$\textit{\ then the state }$y$\textit{\
is more regular, as in} (\ref{7-1-0}).

\medskip

\noindent \textbf{Proof. }Under the specified hypotheses, problem (\ref{1})-(%
\ref{3}) has a unique nonnegative solution given by Theorem 2.4\textbf{. }%
Then, $J(u)\geq 0$, its infimum exists and it is nonnegative. Let us denote
it by $d.$

For not overloading the notations we shall drop the superscript $^{u}.$

We consider a minimizing sequence $(u_{n})_{n\geq 1},$ $u_{n}\in U$ which
satisfies
\begin{eqnarray}
d &\leq &\frac{\lambda _{1}}{2}\left( \int_{Q}u_{n}y_{nx}dxdt-M_{f}\right)
^{2}+\frac{\lambda _{2}}{2}\left( \int_{0}^{L}y_{n}(T,x)dx-M_{T}\right) ^{2}
\label{10} \\
&&+\frac{\lambda _{3}}{2}\left( \int_{Q}y_{n}dxdt-M\right) ^{2}\leq d+\frac{1%
}{n},  \notag
\end{eqnarray}%
where the corresponding state $y_{n}$ is the solution to (\ref{1})-(\ref{3})
(equivalently to (\ref{18-03})) with $u$ replaced by $u_{n}$. By Theorem
2.4,\ it exists for each $n$, it is unique and satisfies%
\begin{equation}
\sup_{t\in \lbrack 0,T]}\left\Vert y_{n}(t)\right\Vert
_{H}^{2}+\int_{0}^{T}\left\Vert y_{n}(t)\right\Vert _{V_{u_{n}}}^{2}dt\leq C,%
\mbox{ for any }t\in \lbrack 0,T],  \label{11}
\end{equation}%
with $C$ a positive constant independent of $n,$ by (\ref{7-2}).

Since $u_{n}\in U,$ we deduce that there exist subsequences (denoted still
by the subscript $_{n})$ such that
\begin{equation*}
u_{n}\rightarrow u\mbox{ weak* in }L^{\infty }(0,L),\mbox{ as }n\rightarrow
\infty ,
\end{equation*}%
\begin{equation*}
u_{nx}\rightarrow u_{x}\mbox{ weak* in }L^{\infty }(0,L),\mbox{ as }%
n\rightarrow \infty ,
\end{equation*}%
\begin{equation*}
y_{n}\rightarrow y\mbox{ weak* in }L^{\infty }(0,T;L^{2}(0,L)),\mbox{ as }%
n\rightarrow \infty ,
\end{equation*}%
\begin{equation*}
y_{n}(T)\rightarrow \zeta \mbox{ weakly in }L^{2}(0,L),\mbox{ as }%
n\rightarrow \infty .
\end{equation*}%
By (\ref{11}), $(\sqrt{u_{n}}y_{nx})_{n}$ is bounded in $L^{2}(Q)$ and there
exists $\xi \in L^{2}(Q)$ such that, on a subsequence (denoted still by the
subscript $_{n})$%
\begin{equation*}
\xi _{n}:=\sqrt{u_{n}}y_{nx}\rightarrow \xi \mbox{ weakly in }L^{2}(Q),\mbox{
as }n\rightarrow \infty .
\end{equation*}

Let $\delta $ be positive, arbitrary. Then $\xi _{n}\rightarrow \xi $ weakly
in $L^{2}(0,T;L^{2}(\Omega _{\delta }))$ too$,$ where $\Omega _{\delta
}=(0,x_{0}-\delta )\cup (x_{0}+\delta ,L).$

The previous first two convergences for $u_{n}$ imply that
\begin{equation*}
u_{n}\rightarrow u\mbox{ uniformly on }[0,L],\mbox{ as }n\rightarrow \infty ,
\end{equation*}%
and therefore the sequence $y_{nx}=\frac{\xi _{n}}{\sqrt{u_{n}}}$ is bounded
in $L^{2}(0,T;L^{2}(\Omega _{\delta }))$ and%
\begin{equation*}
y_{nx}=\frac{\xi _{n}}{\sqrt{u_{n}}}\rightarrow \frac{\xi }{\sqrt{u}}\mbox{
weakly in }L^{2}(0,T;L^{2}(\Omega _{\delta })),\mbox{ as }n\rightarrow
\infty .
\end{equation*}%
On the other hand, $y_{nx}\rightarrow y_{x}$ in the sense of distributions.
We conclude that $\xi =\sqrt{u}y_{x}$ a.e. on $(0,T)\times \Omega _{\delta }$
and since $\delta $ is arbitrary we finally get $\xi =\sqrt{u}y_{x}$ a.e. on
$Q.$ It follows that $y\in L^{2}(0,T;H_{u}^{1}(0,L))$ and so
\begin{equation*}
y_{n}\rightarrow y\mbox{ weakly in }L^{2}(0,T;H_{u}^{1}(0,L)),\mbox{ as }%
n\rightarrow \infty .
\end{equation*}%
Next, by the definition of $A,$ see (\ref{18-00}), we have that%
\begin{equation*}
Ay_{n}\rightarrow Ay\mbox{ weakly in }L^{2}(0,T;(H_{u}^{1}(0,L))^{\prime }),%
\mbox{ as }n\rightarrow \infty ,
\end{equation*}%
and by (\ref{18-03}) we deduce that $\left( \frac{dy_{n}}{dt}\right) _{n}$
is bounded in $L^{2}(0,T;(H_{u}^{1}(0,L))^{\prime }),$ so that by selecting
a subsequence we have%
\begin{equation*}
\frac{dy_{n}}{dt}\rightarrow \frac{dy}{dt}\mbox{ weakly in }%
L^{2}(0,T;(H_{u}^{1}(0,L))^{\prime }),\mbox{ as }n\rightarrow \infty .
\end{equation*}%
Moreover, since $y_{n}$ is the solution to (\ref{18-03}) we get that
\begin{eqnarray}
y_{n}(t) &=&y_{0}+\int_{0}^{t}(-Ay_{n}(s)+f(s))ds  \label{si-200} \\
&\rightarrow &y_{0}+\int_{0}^{t}(-Ay(s)+f(s))ds,\mbox{ weakly in }%
(H_{u}^{1}(0,L))^{\prime },\mbox{ for all }t\in \lbrack 0,T].  \notag
\end{eqnarray}%
From here we get that%
\begin{equation}
\int_{0}^{L}y_{n}(t)\phi _{0}dx=\int_{0}^{L}y_{0}\phi
_{0}dx+\int_{0}^{t}\left\langle -Ay_{n}(s)+f(s),\phi _{0}\right\rangle
_{V_{u}^{\prime },V_{u}}ds  \label{si-200-0}
\end{equation}%
for any $\phi _{0}\in H_{u}^{1}(0,L)$ and $t\in \lbrack 0,T].$ Passing to
the limit we obtain
\begin{equation*}
l(t)=\lim_{n\rightarrow \infty }\int_{0}^{L}y_{n}(t)\phi
_{0}dx=\int_{0}^{L}y_{0}\phi _{0}dx+\int_{0}^{t}\left\langle
-Ay(s)+f(s),\phi _{0}\right\rangle _{V_{u}^{\prime },V_{u}}ds.
\end{equation*}%
We multiply this relation by $\phi _{1}\in L^{2}(0,T)$ and integrate over $%
(0,T),$ getting%
\begin{equation}
\int_{0}^{T}\phi _{1}(t)l(t)dt=\int_{0}^{T}\left( \int_{0}^{L}y_{0}\phi
_{0}dx+\int_{0}^{t}\left\langle -Ay(s)+f(s),\phi _{0}\right\rangle
_{V_{u}^{\prime },V_{u}}ds\right) \phi _{1}(t)dt.  \label{si-201}
\end{equation}%
We multiply (\ref{si-200}) by $\phi _{1}(t)\phi _{0}(x)$ and integrate over $%
(0,T)\times (0,L),$%
\begin{equation}
\int_{Q}\phi _{0}\phi _{1}y_{n}dxdt=\int_{0}^{T}\left( \int_{0}^{L}y_{0}\phi
_{0}dx+\int_{0}^{t}\left\langle -Ay_{n}(s)+f(s),\phi _{0}\right\rangle
_{V_{u}^{\prime },V_{u}}ds\right) \phi _{1}(t)dt.  \label{si-202}
\end{equation}%
By the weak convergence $y_{n}\rightarrow y$ in $L^{2}(Q)$ we get that
\begin{equation}
\int_{Q}\phi _{0}\phi _{1}ydxdt=\int_{0}^{T}\left( \int_{0}^{L}y_{0}\phi
_{0}dx+\int_{0}^{t}\left\langle -Ay(s)+f(s),\phi _{0}\right\rangle
_{V_{u}^{\prime },V_{u}}ds\right) \phi _{1}(t)dt.  \label{si-203}
\end{equation}%
Comparing (\ref{si-201}) and (\ref{si-203}) we deduce that
\begin{equation*}
\int_{0}^{T}\phi _{1}(t)l(t)dt=\int_{Q}\phi _{0}\phi _{1}ydxdt,\mbox{ for
any }\phi _{0}\in H_{u}^{1}(0,L),\phi _{1}\in L^{2}(0,T),
\end{equation*}%
hence
\begin{equation*}
l(t)=\lim_{n\rightarrow \infty }\int_{0}^{L}y_{n}(t)\phi
_{0}dx=\int_{0}^{L}y(t)\phi _{0}dx,\mbox{ for any }\phi _{0}\in
H_{u}^{1}(0,L),\mbox{ }t\in \lbrack 0,T].
\end{equation*}%
Therefore,
\begin{equation*}
y_{0}=y_{n}(0)\rightarrow y(0)\mbox{, }y_{n}(T)\rightarrow y(T)\mbox{ weakly
in }(H_{u}^{1}(0,L))^{\prime },\mbox{ as }n\rightarrow \infty
\end{equation*}%
and so by the limit uniqueness $\zeta =y(T)$ and $y(0)=y_{0}$ a.e. on $(0,L)$%
.

Now, $y_{n}$ satisfies (\ref{18-11})
\begin{equation*}
\int_{0}^{T}\left\langle \frac{dy_{n}}{dt}(t),\psi (t)\right\rangle
_{V_{u}^{\prime },V_{u}}dt+\int_{Q}u_{n}y_{nx}\psi _{x}dxdt=\int_{Q}f\psi
dxdt,
\end{equation*}%
for any $\psi \in L^{2}(0,T;H_{u}^{1}(0,L))$ and passing to the limit as $%
n\rightarrow \infty $ we get that $y$ satisfies (\ref{18-11}), too. All
these assertions prove that $y$ is the solution to (\ref{1})-(\ref{3})
corresponding to $u.$

Next, if $y_{0}\in H_{u}^{1}(0,L),$ the solution $y$ previously obtained has
the regularity (\ref{7-1-0}) according to Theorem 2.4$.$

Finally, we pass to the limit in (\ref{10}) as $n\rightarrow \infty ,$ on
the basis of the weakly lower semicontinuity of each convex term in $%
J(u_{n}),$ and get that $d=J(u).$

Since $U$ is closed, then $u(x)\in \lbrack u_{m}(x),u_{M}(x)]$ which implies
by (\ref{3-1}) and (\ref{3-2}) that $u(x)>0$ on $[0,L]\backslash \{x_{0}\},$
$u(x_{0})=0$ and $\frac{1}{u}\notin L^{1}(0,L),$ so that $u\in U$ and the
corresponding operator $A_{1}y=(uy_{x})_{x}$ is strongly degenerate. \hfill $%
\square $

\subsection{Optimality conditions in the homogeneous Dirichlet case}

\noindent \textbf{Proposition 2.6. }\textit{Let }$(u^{\ast },y^{\ast })$%
\textit{\ be a solution to }$(P).$\textit{\ Then }$u^{\ast }$ \textit{%
satisfies} \textit{the necessary condition }%
\begin{equation}
\int_{Q}(u^{\ast }-u)y_{x}^{\ast }\left[ p_{x}+\lambda _{1}\left(
\int_{Q}u^{\ast }y_{x}^{\ast }dxdt-M_{f}\right) \right] dxdt\leq 0
\label{25}
\end{equation}%
\textit{for all} $u\in U,$ \textit{where }$p$\textit{\ is the solution to}
\begin{eqnarray}
&&\frac{\partial p}{\partial t}+(u^{\ast }p_{x})_{x}  \label{19} \\
&=&-\lambda _{1}\left( \int_{Q}u^{\ast }y_{x}^{\ast }dxdt-M_{f}\right)
+\lambda _{3}\left( \int_{Q}y^{\ast }dxdt-M\right) \mbox{ \textit{in} }Q,
\notag
\end{eqnarray}%
\begin{equation}
p(T,x)=-\lambda _{2}\left( \int_{0}^{L}y^{\ast }(T,x)dx-M_{T}\right) \mbox{
\textit{in} }(0,L),  \label{20}
\end{equation}%
\begin{equation}
p(t,0)=p(t,L)=0\mbox{ \textit{in} }(0,T).  \label{21}
\end{equation}

\medskip

\noindent \textbf{Proof. }Let $(u^{\ast },y^{\ast })$ be a solution to $(P)$%
, $\lambda >0,$ $u\in U$ and denote
\begin{equation*}
u^{\lambda }(x)=u^{\ast }(x)+\lambda v(x),
\end{equation*}%
where
\begin{equation}
v(x)=u(x)-u^{\ast }(x),\mbox{ }u\in U.  \label{15}
\end{equation}%
It is obvious that $v\in W^{1,\infty }(0,L),$ $v(x_{0})=0,$ $%
u_{m}(x)-u_{M}(x)\leq v(x)\leq u(x),$ $v(0)=v(L)=0$ and $\frac{1}{v}\notin
L^{1}(0,L).$ We introduce the system
\begin{equation}
\frac{\partial Y}{\partial t}-(u^{\ast }Y_{x})_{x}=(vy_{x}^{\ast })_{x}\mbox{
in }Q,  \label{16}
\end{equation}%
\begin{equation}
Y(0,x)=0\mbox{ in }(0,L),  \label{17}
\end{equation}%
\begin{equation}
Y(t,0)=Y(t,L)=0\mbox{ in }(0,T).  \label{18}
\end{equation}%
We note that since $y^{\ast }\in L^{2}(0,T;H_{u^{\ast }}^{1}(0,L)),$ after a
few calculations it follows that%
\begin{equation*}
\int_{Q}v^{2}(y_{x}^{\ast })^{2}dxdt\leq \left\Vert u_{M}\right\Vert
_{C[0,L]}\left\Vert \alpha \right\Vert _{L^{\infty
}(0,L)}\int_{0}^{T}\left\Vert \sqrt{u^{\ast }}y_{x}^{\ast }(t)\right\Vert
_{H}^{2}dt.
\end{equation*}%
Hence, $vy_{x}^{\ast }\in L^{2}(0,T;L^{2}(0,L))$ and $g=(vy_{x}^{\ast
})_{x}\in L^{2}(0,T;(H_{u^{\ast }}^{1}(0,L))^{\prime }),$ since for any $%
\psi \in H_{u^{\ast }}^{1}(0,L),$ and a.e. $t\in (0,T)$ we have%
\begin{equation*}
\left\vert \left\langle (vy_{x}^{\ast })_{x}(t),\psi \right\rangle
_{V_{u^{\ast }}^{\prime },V_{u^{\ast }}}\right\vert =\left\vert
-\int_{0}^{L}vy_{x}^{\ast }(t)\psi _{x}dx\right\vert \leq \left\Vert \alpha
\right\Vert _{L^{\infty }(0,L)}\left\Vert y^{\ast }(t)\right\Vert
_{V_{u^{\ast }}}\left\Vert \psi \right\Vert _{V_{u^{\ast }}}.
\end{equation*}

We state that (\ref{16})-(\ref{18}) has a unique solution
\begin{equation}
Y\in C([0,T];L^{2}(0,L))\cap L^{2}(0,T;V_{u^{\ast }})\cap
W^{1,2}([0,T];V_{u^{\ast }}^{\prime }).  \label{18-0}
\end{equation}%
This follows as in the first part of Theorem 2.4 by defining the operator $A$
from $H_{u^{\ast }}^{1}(0,L)$ to $(H_{u^{\ast }}^{1}(0,L))^{\prime },$ by%
\begin{equation}
\left\langle Az,\psi \right\rangle _{V_{u^{\ast }}^{\prime },V_{u^{\ast
}}}=\int_{0}^{L}u^{\ast }z_{x}\psi _{x}dx,\mbox{ for any }\psi \in
V_{u^{\ast }}.  \label{18-1}
\end{equation}

Moreover, denoting by $y^{\lambda }(t,x)$ the solution to (\ref{1})-(\ref{3}%
) corresponding to $u^{\lambda }(x),$ one can prove that actually
\begin{equation*}
Y(t,x)=\lim_{\lambda \rightarrow 0}\frac{y^{\lambda }(t,x)-y^{\ast }(t,x)}{%
\lambda },
\end{equation*}%
so that (\ref{16})-(\ref{18}) is the system of first order variations.

We introduce the dual system (\ref{19})-(\ref{21})$.$ This system has a
unique solution
\begin{multline*}
p\in C([0,T];H_{u^{\ast }}^{1}(0,L))\cap L^{2}(0,T;H_{u^{\ast
}}^{2}(0,L))\cap W^{1,2}([0,T];L^{2}(0,L)),\mbox{ } \\
(u^{\ast }p_{x})_{x}\in L^{2}(Q),
\end{multline*}%
given still by Theorem 2.4, second part (after making the transformation $%
t^{\prime }=T-t).$

Now, we write that $(u^{\ast },y^{\ast })$ is a solution to $(P),$ that is
\begin{equation*}
J(u^{\ast })\leq J(u),\mbox{ for all }u\in U,
\end{equation*}%
and, in particular, for $u=u^{\lambda }.$ After some algebra we get%
\begin{equation}
\lambda _{1}I_{f}\int_{Q}(-u_{x}^{\ast }Y+vy_{x}^{\ast })dxdt+\lambda
_{2}I_{T}\int_{0}^{L}Y(T,x)dx+\lambda _{3}I_{M}\int_{Q}Ydxdt\geq 0,
\label{22}
\end{equation}%
where
\begin{equation}
I_{f}=\int_{Q}u^{\ast }y_{x}^{\ast }dxdt-M_{f},\mbox{ }I_{T}=\int_{0}^{L}y^{%
\ast }(T,x)dx-M_{T},\mbox{ }I_{M}=\int_{Q}y^{\ast }dxdt-M.  \label{22-0}
\end{equation}%
We test (\ref{16}) by $p(t)$ and integrate over $(0,T).$ After some
calculations we get
\begin{eqnarray*}
&&\int_{Q}-\left( p_{t}+(u^{\ast }(x)p_{x})_{x}\right)
Ydxdt+\int_{0}^{L}p(T,x)Y(T,x)dx \\
&=&\int_{0}^{T}\left\langle (vy_{x}^{\ast }(t))_{x},p(t)\right\rangle
_{V_{u^{\ast }}^{\prime },V_{u^{\ast }}}dt,
\end{eqnarray*}%
which yields, by (\ref{19})-(\ref{21}),
\begin{equation}
\int_{Q}(-\lambda _{1}I_{f}u_{x}^{\ast }+\lambda _{3}I_{M})Ydxdt+\lambda
_{2}I_{T}\int_{0}^{L}Y(T,x)dx=\int_{Q}vy_{x}^{\ast }p_{x}dxdt.  \label{23}
\end{equation}%
Comparing with (\ref{22}) it follows that%
\begin{equation*}
\int_{Q}vy_{x}^{\ast }p_{x}dxdt+\int_{Q}\lambda _{1}I_{f}vy_{x}^{\ast
}dxdt\geq 0
\end{equation*}%
with $v=u-u^{\ast },$ for all $u\in U,$ and this implies (\ref{25}), as
claimed.\hfill $\square $

\subsection{Approximating problem}

In order to give a better characterization of the optimality condition (\ref%
{25}) we determine an approximating form of it. To this end, for $%
\varepsilon >0,$ we introduce an approximating problem $(P_{\varepsilon })$
involving a nondegenerate state equation. The approximating optimality
condition may be written more explicitly due to the better regularity of the
approximating state and dual variable. Then we show that $(P_{\varepsilon })$
tends in some sense to $(P).$ Namely, we show that a sequence of solutions
to $(P_{\varepsilon })$ tends to a solution to $(P),$ as $\varepsilon
\rightarrow 0$.

We introduce
\begin{equation}
\mbox{Minimize }J(u)\mbox{ for all }u\in U_{\varepsilon },
\tag{$P_{\varepsilon }$}
\end{equation}%
subject to the state system (\ref{1})-(\ref{3}), where%
\begin{eqnarray}
U_{\varepsilon } &=&\{u\in W^{1,\infty }(0,L);\mbox{ }u_{m}(x)+\varepsilon
\leq u(x)\leq u_{M}(x)+2\varepsilon ,\mbox{ }  \label{Ueps} \\
u(0) &=&u_{0}^{\varepsilon },\mbox{ }u(L)=u_{L}^{\varepsilon },\mbox{ }%
\left\vert u_{x}(x)\right\vert \leq u_{\infty }\mbox{ a.e. }x\in (0,L)\}.
\notag
\end{eqnarray}%
All hypotheses made in (\ref{3-1})-(\ref{3-2-0}) remain the same and we note
that if $u\in U_{\varepsilon },$ then%
\begin{equation*}
\varepsilon \leq u(x_{0})\leq 2\varepsilon ,
\end{equation*}%
\begin{equation*}
u_{m}(0)+\varepsilon \leq u_{0}^{\varepsilon }\leq u_{M}(0)+2\varepsilon ,%
\mbox{ \ }u_{m}(L)+\varepsilon \leq u_{L}^{\varepsilon }\leq
u_{M}(L)+2\varepsilon .
\end{equation*}

For all $u\in U_{\varepsilon },$ $u(x)\geq u_{m}(x)+\varepsilon \geq
\varepsilon ,$ and then system (\ref{1})-(\ref{3}) with $u\in U_{\varepsilon
}$ is nondegenerate.

By the general results concerning nondegenerate evolution equations in
Hilbert spaces, if $y_{0}\in L^{2}(0,L),$ and $f\in L^{2}(Q),$ problem (\ref%
{1})-(\ref{3}) with $u(x)\geq \varepsilon $ has a unique solution
\begin{equation}
y_{\varepsilon }\in C([0,T];L^{2}(0,L))\cap L^{2}(0,T;H_{0}^{1}(0,L))\cap
W^{1,2}([0,T];H^{-1}(0,L)),  \label{y-nondeg}
\end{equation}%
(see \cite{lions}, p. 163) and $y_{\varepsilon }$ satisfies the estimate%
\begin{equation}
\sup_{t\in \lbrack 0,T]}\left\Vert y_{\varepsilon }(t)\right\Vert
_{H}^{2}+\int_{0}^{T}\left\Vert \sqrt{u}y_{\varepsilon x}(t)\right\Vert
_{H}^{2}dt\leq C,  \label{31}
\end{equation}%
with $C$ a positive constant depending on the data and independent of $%
\varepsilon .$

We denote by $y_{\varepsilon x}$ the derivative of $y_{\varepsilon }$ with
respect to $x.$

Obviously, the control problem $(P_{\varepsilon })$ has at least a solution $%
(u_{\varepsilon },y_{\varepsilon }),$ with $u_{\varepsilon }\in
U_{\varepsilon }$ and $y_{\varepsilon }$ (corresponding to $u_{\varepsilon
}) $ satisfying (\ref{31}). We prove the convergence result $(P_{\varepsilon
})\rightarrow (P)$ as $\varepsilon \rightarrow 0,$ on the basis of the
following lemma.

\medskip

\noindent \textbf{Lemma 2.7. }\textit{Let }$y_{0}\in L^{2}(0,L),$ $f\in
L^{2}(Q)$, $y_{0}\geq 0$ \textit{a.e. on }$(0,L),$\textit{\ }$f\geq 0$%
\textit{\ a.e. on }$Q.$\textit{\ Let} $(u_{\varepsilon },y_{\varepsilon })$%
\textit{\ be a sequence of solutions to }$(P_{\varepsilon })$\textit{\ such
that }%
\begin{equation}
u_{\varepsilon }\rightarrow u\mbox{ \textit{uniformly in} }[0,L],\mbox{
\textit{as} }\varepsilon \rightarrow 0,  \label{32}
\end{equation}%
\begin{equation}
u_{\varepsilon x}\rightarrow u_{x}\mbox{ \textit{weak* in} }L^{\infty }(0,L),%
\mbox{ \textit{as} }\varepsilon \rightarrow 0.  \label{33}
\end{equation}%
\textit{Then, }$u\in U,$%
\begin{equation}
y_{\varepsilon }\rightarrow y\mbox{ \textit{weakly in} }%
L^{2}(0,T;H_{u}^{1}(0,L))\cap W^{1,2}([0,T];(H_{u}^{1}(0,L))^{\prime }),%
\mbox{ \textit{as} }\varepsilon \rightarrow 0,  \label{34}
\end{equation}%
\begin{equation}
y_{\varepsilon }(T)\rightarrow y(T)\mbox{ \textit{weakly} \textit{in} }%
L^{2}(0,L),\mbox{ \textit{as} }\varepsilon \rightarrow 0,  \label{34-2}
\end{equation}%
\textit{and }$y$\textit{\ is the solution to} (\ref{1})-(\ref{3}) \textit{%
corresponding to }$u.$\textit{\ }

\medskip

\noindent \textbf{Proof. }It is easily seen that\textbf{\ }$%
\lim\limits_{\varepsilon \rightarrow 0}u_{\varepsilon }=u\in U$ defined by (%
\ref{U}). Then, $y_{\varepsilon }$ is the solution to the nondegenerate
problem (\ref{1})-(\ref{3}) corresponding to $u_{\varepsilon }$ and $%
y_{\varepsilon }$ satisfies (\ref{31}). We get (on subsequences denoted
still by the subscript $_{\varepsilon })$ that
\begin{equation*}
y_{\varepsilon }\rightarrow y\mbox{ weakly in }L^{2}(0,T;L^{2}(0,L)),\mbox{
as }\varepsilon \rightarrow 0,
\end{equation*}%
\begin{equation*}
\sqrt{u_{\varepsilon }}y_{\varepsilon }\rightarrow \xi \mbox{ weakly in }%
L^{2}(0,T;L^{2}(0,L)),\mbox{ as }\varepsilon \rightarrow 0,
\end{equation*}%
\begin{equation*}
y_{\varepsilon }(T)\rightarrow \zeta \mbox{ weakly in }L^{2}(0,L),\mbox{ as }%
\varepsilon \rightarrow 0.
\end{equation*}%
We continue the proof in a similar way as in Theorem 2.5 and we get that $%
\xi =\sqrt{u}y$ a.e. on $Q$ and so $y\in L^{2}(0,T;H_{u}^{1}(0,L)).$ By (\ref%
{18-00}) we deduce that
\begin{equation*}
Ay_{\varepsilon }\rightarrow Ay\mbox{ weakly in }%
L^{2}(0,T;(H_{u}^{1}(0,L))^{\prime }),\mbox{ as }\varepsilon \rightarrow 0
\end{equation*}%
and by (\ref{18-03})
\begin{equation*}
\frac{dy_{\varepsilon }}{dt}\rightarrow \frac{dy}{dt}\mbox{ weakly in }%
L^{2}(0,T;(H_{u}^{1}(0,L))^{\prime }),\mbox{ as }\varepsilon \rightarrow 0.
\end{equation*}%
Following the arguments in Theorem 2.5 we prove (\ref{34-2}) and that $y$ is
the solution to (\ref{1})-(\ref{3}) corresponding to $u.$\hfill $\square $

\medskip

\noindent \textbf{Theorem 2.8. }\textit{Let }$y_{0}\in L^{2}(0,L)$, $f\in
L^{2}(Q),$ $y_{0}\geq 0$ \textit{a.e. on }$(0,L)$\textit{, }$f\geq 0$\textit{%
\ a.e. on} $Q.$ \textit{Let }$(u_{\varepsilon }^{\ast },y_{\varepsilon
}^{\ast })_{\varepsilon >0}$\textit{\ be a sequence of solutions to }$%
(P_{\varepsilon }).$\textit{\ Then }(\textit{on subsequences denoted still
by }$_{\varepsilon }$)\textit{\ we have }%
\begin{equation}
u_{\varepsilon }^{\ast }\rightarrow u^{\ast }\mbox{ \textit{uniformly in} }%
[0,L],\mbox{ \textit{as} }\varepsilon \rightarrow 0,  \label{32-0}
\end{equation}%
\begin{equation}
u_{\varepsilon x}^{\ast }\rightarrow u_{x}^{\ast }\mbox{ \textit{weak* in} }%
L^{\infty }(0,L),\mbox{ \textit{as} }\varepsilon \rightarrow 0,  \label{33-0}
\end{equation}%
\begin{equation}
y_{\varepsilon }^{\ast }\rightarrow y^{\ast }\mbox{ \textit{weakly in} }%
L^{2}(0,T;H_{u^{\ast }}^{1}(0,L))\cap W^{1,2}([0,T];(H_{u^{\ast
}}^{1}(0,L))^{\prime }),\mbox{ \textit{as} }\varepsilon \rightarrow 0,
\label{34-0}
\end{equation}%
\begin{equation}
y_{\varepsilon }(T)\rightarrow y(T)\mbox{ \textit{weakly in} }L^{2}(0,L),%
\mbox{ \textit{as} }\varepsilon \rightarrow 0.  \label{34-000}
\end{equation}%
\textit{Moreover, }$y^{\ast }$\textit{\ is the solution to} (\ref{1})-(\ref%
{3}) \textit{corresponding to }$u^{\ast }$\textit{\ and }$(u^{\ast },y^{\ast
})$\textit{\ is a solution to }$(P).$\textit{\ }

\medskip

\noindent \textbf{Proof. }Let $(u_{\varepsilon }^{\ast },y_{\varepsilon
}^{\ast })$ be a solution to $(P_{\varepsilon })$, i.e.,%
\begin{equation*}
J(u_{\varepsilon }^{\ast })\leq J(u_{\varepsilon }),\mbox{ for all }%
u_{\varepsilon }\in U_{\varepsilon }.
\end{equation*}%
Under the hypotheses for\textbf{\ }$y_{0}$ it follows that (\ref{1})-(\ref{3}%
) with $u_{\varepsilon }\in U_{\varepsilon }$ has a unique solution $%
y_{\varepsilon }\in C([0,T];L^{2}(0,L))\cap W^{1,2}([0,T];(H_{u_{\varepsilon
}}^{1}(0,L))^{\prime })\cap L^{2}(0,T;H_{u_{\varepsilon }}^{1}(0,L))$ and
\begin{eqnarray}
&&\frac{\lambda _{1}}{2}\left( \int_{Q}u_{\varepsilon }^{\ast
}y_{\varepsilon x}^{\ast }dxdt-M_{f}\right) ^{2}+\frac{\lambda _{2}}{2}%
\left( \int_{0}^{L}y_{\varepsilon }^{\ast }(T,x)dx-M_{T}\right) ^{2}
\label{35} \\
&&+\frac{\lambda _{3}}{2}\left( \int_{Q}y_{\varepsilon }^{\ast
}dxdt-M\right) ^{2}  \notag \\
&\leq &\frac{\lambda _{1}}{2}\left( \int_{Q}u_{\varepsilon }y_{\varepsilon
x}dxdt-M_{f}\right) ^{2}+\frac{\lambda _{2}}{2}\left(
\int_{0}^{L}y_{\varepsilon }(T,x)dx-M_{T}\right) ^{2}  \notag \\
&&+\frac{\lambda _{3}}{2}\left( \int_{Q}y_{\varepsilon }dxdt-M\right) ^{2}
\notag
\end{eqnarray}%
for all $u_{\varepsilon }\in U_{\varepsilon }.$ Relations (\ref{32-0})-(\ref%
{33-0}) follow from $u_{\varepsilon }^{\ast }\in U_{\varepsilon }$ and so
the sequence of the corresponding states ($y_{\varepsilon }^{\ast
})_{\varepsilon }$ converges on a subsequence to $y^{\ast }$ the solution to
(\ref{1})-(\ref{3}) corresponding to $u^{\ast },$ as established by Lemma
2.7. In particular, we have (\ref{34}) which implies that
\begin{equation*}
u_{\varepsilon }y_{\varepsilon x}\rightarrow u^{\ast }y_{x}^{\ast }\mbox{
weakly in }L^{2}(Q),\mbox{ as }\varepsilon \rightarrow 0.
\end{equation*}%
Similarly, $u_{\varepsilon }\in U_{\varepsilon }$ implies that $%
u_{\varepsilon }\rightarrow u\in U$ uniformly on $[0,L]$ as $\varepsilon
\rightarrow 0$, and by Lemma 2.7, (\ref{34})-(\ref{34-2}), we get that $%
(y_{\varepsilon })_{\varepsilon }$ is convergent to $y$ which is the
solution to (\ref{1})-(\ref{3}) corresponding to $u.$

Passing to the limit in (\ref{35}) we get
\begin{multline*}
\frac{\lambda _{1}}{2}\left( \int_{Q}u^{\ast }y_{x}^{\ast }dxdt-M_{f}\right)
^{2}+\frac{\lambda _{2}}{2}\left( \int_{0}^{L}y^{\ast }(T,x)dx-M_{T}\right)
^{2}+\frac{\lambda _{3}}{2}\left( \int_{Q}y^{\ast }dxdt-M\right) ^{2} \\
\leq \liminf\limits_{\varepsilon \rightarrow 0}\left( \frac{\lambda _{1}}{2}%
\left( \int_{Q}u_{\varepsilon }^{\ast }y_{\varepsilon x}^{\ast
}dxdt-M_{f}\right) ^{2}+\frac{\lambda _{2}}{2}\left(
\int_{0}^{L}y_{\varepsilon }^{\ast }(T,x)dx-M_{T}\right) ^{2}\right. \\
\left. +\frac{\lambda _{3}}{2}\left( \int_{Q}y_{\varepsilon }^{\ast
}dxdt-M\right) ^{2}\right) \\
\leq \limsup\limits_{\varepsilon \rightarrow 0}\left( \frac{\lambda _{1}}{2}%
\left( \int_{Q}u_{\varepsilon }y_{\varepsilon x}dxdt-M_{f}\right) ^{2}+\frac{%
\lambda _{2}}{2}\left( \int_{0}^{L}y_{\varepsilon }(T,x)dx-M_{T}\right)
^{2}\right. \\
\left. +\frac{\lambda _{3}}{2}\left( \int_{Q}y_{\varepsilon }dxdt-M\right)
^{2}\right) \\
\leq \frac{\lambda _{1}}{2}\left( \int_{Q}uy_{x}dxdt-M_{f}\right) ^{2}+\frac{%
\lambda _{2}}{2}\left( \int_{0}^{L}y(T,x)dx-M_{T}\right) ^{2}+\frac{\lambda
_{3}}{2}\left( \int_{Q}ydxdt-M\right) ^{2},
\end{multline*}%
for all $u\in U.$ This implies that $(u^{\ast },y^{\ast })$ is a solution to
$(P).$\hfill $\square $

\subsection{Approximating optimality conditions}

Let $K$ be a closed convex subset of a Banach space $X$ having the dual $%
X^{\prime }.$ We recall (see \cite{vb3}, p. 4-5) that the indicator function
of $K$ is
\begin{equation*}
I_{K}(\xi )=\left\{
\begin{array}{l}
0,\mbox{ \ \ \ \ if }\xi \in K, \\
+\infty ,\mbox{ if }\xi \notin K,%
\end{array}%
\right.
\end{equation*}%
and the subdifferential of $I_{K}$ coincides with the normal cone to $K$ at $%
\xi ,$
\begin{equation*}
\partial I_{K}(\xi )=N_{K}(\xi )=\{\xi ^{\ast }\in X^{\prime };\mbox{ }%
\left\langle \xi ^{\ast },\xi \right\rangle _{X^{\prime },X}\geq 0\}.
\end{equation*}%
For $\theta :\mathbb{R}\rightarrow \lbrack -1,1],$ let us denote by sign$%
~\theta $ the graph
\begin{equation*}
\mbox{sign}~\theta =\left\{
\begin{array}{l}
1,\mbox{ \ \ \ \ \ \ \ on }\{x;\mbox{ }\theta (x)>0\}, \\
\lbrack -1,1],\mbox{ on }\{x;\mbox{ }\theta (x)=0\}, \\
-1,\mbox{ \ \ \ \ \ on }\{x;\mbox{ }\theta (x)<0\}.%
\end{array}%
\right.
\end{equation*}

\medskip

\noindent \textbf{Proposition 2.9. }\textit{Let us assume the hypotheses as
in }Theorem 2.8.\textit{\ Then, the approximating optimality condition reads
as}%
\begin{equation}
\int_{0}^{L}(u_{\varepsilon }^{\ast }-u_{\varepsilon })(x)\Phi (x)dx\geq 0,%
\mbox{ \textit{for all} }u_{\varepsilon }\in U_{\varepsilon },  \label{100}
\end{equation}%
\textit{\ where}
\begin{equation}
\Phi (x):=-\int_{0}^{T}y_{\varepsilon x}^{\ast }(t,x)\left\{ p_{\varepsilon
x}(t,x)+\lambda _{1}\left( \int_{Q}u_{\varepsilon }^{\ast }y_{\varepsilon
x}^{\ast }dxdt-M_{f}\right) \right\} dt.  \label{41}
\end{equation}%
\textit{Moreover, }$\Phi $ \textit{has the representation}
\begin{equation}
\Phi (x)=-\rho ^{\prime }(x)+\mu (x)\mbox{ \textit{a.e. }}x\in (0,L),
\label{71}
\end{equation}%
\begin{equation}
u_{\varepsilon x}^{\ast }(x)\in u_{\infty }\mbox{sign(}\rho (x))\mbox{
\textit{a.e. }}x\in (0,L),  \label{72}
\end{equation}%
\textit{where }$\mu ,$ $\rho \in L^{1}(0,L),$
\begin{equation}
\left\{
\begin{array}{l}
\mu (x)\leq 0\mbox{ \ \textit{a.e. in} }\{x\in (0,L);\mbox{ }u_{\varepsilon
}^{\ast }(x)=u_{m}(x)+\varepsilon \} \\
\mu (x)=0\mbox{ \ \textit{a.e. in} }\{x\in (0,L);\mbox{ }u_{\varepsilon
}^{\ast }(x)\in (u_{m}(x)+\varepsilon ,u_{M}(x)+2\varepsilon )\} \\
\mu (x)\geq 0\mbox{ \ \textit{a.e. in} }\{x\in (0,L);\mbox{ }u_{\varepsilon
}^{\ast }(x)=u_{M}(x)+2\varepsilon \},%
\end{array}%
\right.  \label{70}
\end{equation}%
\textit{\ and }%
\begin{equation}
\rho (x)=\left\{
\begin{array}{l}
0\mbox{ \ \ \ \ \ \ \ \ \ \ \ \ \textit{a.e. in} }\{x\in (0,L);\left\vert
u_{\varepsilon x}^{\ast }(x)\right\vert <u_{\infty }\mbox{ \textit{a.e.}}\}
\\
\nu (x)u_{\varepsilon x}^{\ast }(x)\mbox{ \textit{a.e. in} }\{x\in
(0,L);\left\vert u_{\varepsilon x}^{\ast }(x)\right\vert =u_{\infty }\mbox{
\textit{a.e.}}\}%
\end{array}%
\right.  \label{44}
\end{equation}%
\textit{with }$\nu \in L^{1}(0,L),$\textit{\ }$\nu \geq 0$\textit{\ a.e. }$%
x\in (0,L).$

\medskip

\noindent \textbf{Proof. }The system in variations, the dual system and the
optimality conditions are similarly obtained as those for $(P).$ Namely, if $%
(u_{\varepsilon }^{\ast },y_{\varepsilon }^{\ast })$ is a solution to $%
(P_{\varepsilon })$ we introduce
\begin{eqnarray}
&&\frac{\partial p_{\varepsilon }}{\partial t}+(u_{\varepsilon }^{\ast
}p_{\varepsilon x})_{x}  \label{37} \\
&=&-\lambda _{1}\left( \int_{Q}u_{\varepsilon }^{\ast }y_{\varepsilon
x}^{\ast }dxdt-M_{f}\right) +\lambda _{3}\left( \int_{Q}y_{\varepsilon
}^{\ast }dxdt-M\right) \mbox{ in }Q,,  \notag
\end{eqnarray}%
\begin{equation}
p_{\varepsilon }(T,x)=-\lambda _{2}\left( \int_{0}^{L}y_{\varepsilon }^{\ast
}(T,x)dx-M_{T}\right) \mbox{ in }(0,L),  \label{38}
\end{equation}%
\begin{equation}
p_{\varepsilon }(t,0)=p_{\varepsilon }(t,L)=0\mbox{ in }(0,T),  \label{39}
\end{equation}%
and deduce that $u_{\varepsilon }^{\ast },$ $y_{\varepsilon }^{\ast }$ and
the dual variable $p_{\varepsilon }$ should satisfy a similar relation to (%
\ref{25}), i.e.,
\begin{equation}
\int_{Q}(u_{\varepsilon }^{\ast }-u_{\varepsilon })y_{\varepsilon x}^{\ast }
\left[ p_{\varepsilon x}+\lambda _{1}\left( \int_{Q}u_{\varepsilon }^{\ast
}y_{\varepsilon x}^{\ast }dxdt-M_{f}\right) \right] dxdt\leq 0  \label{40}
\end{equation}%
for all $u_{\varepsilon }\in U_{\varepsilon }.$ The solution $p_{\varepsilon
}$ is regular, since $u_{\varepsilon }^{\ast }\geq \varepsilon ,$ and%
\begin{equation*}
p_{\varepsilon }\in C([0,T];L^{2}(0,L))\cap L^{2}(0,T;H_{0}^{1}(0,L))\cap
W^{1,2}([0,T];H^{-1}(0,L)),
\end{equation*}%
by the same arguments as for $y_{\varepsilon }$ (see (\ref{y-nondeg})).

The supplementary regularity of $y_{\varepsilon }$ and $p_{\varepsilon }$
implies that (\ref{40}) can be still written
\begin{equation}
\int_{0}^{L}(u_{\varepsilon }^{\ast }-u_{\varepsilon })\left(
\int_{0}^{T}y_{\varepsilon x}^{\ast }\left[ p_{\varepsilon x}+\lambda
_{1}\left( \int_{Q}u_{\varepsilon }^{\ast }y_{\varepsilon }^{\ast
}dxdt-M_{f}\right) \right] dt\right) dx\leq 0,  \label{40-0}
\end{equation}%
for all $u_{\varepsilon }\in U_{\varepsilon }.$ We see that $\Phi \in
L^{1}(0,L).$ Then, (\ref{100}) and (\ref{41}) imply that
\begin{equation}
\Phi (x)\in N_{U_{\varepsilon }}(u_{\varepsilon }^{\ast }).  \label{40-1}
\end{equation}%
Here, $N_{U_{\varepsilon }}(u_{\varepsilon }^{\ast })$ means the normal cone
in $L^{1}(0,L)$ to $U_{\varepsilon }$ at $u_{\varepsilon }^{\ast }\in
U_{\varepsilon }$, that is%
\begin{equation*}
N_{U_{\varepsilon }}(u_{\varepsilon }^{\ast })=\{w\in L^{1}(0,L);\mbox{ }%
\int_{0}^{L}(u_{\varepsilon }^{\ast }-u_{\varepsilon })wdx\geq 0,\mbox{ for
all }u_{\varepsilon }\in U_{\varepsilon }\}.
\end{equation*}

We give further a representation of $N_{U_{\varepsilon }}(u_{\varepsilon
}^{\ast }),$ more exactly we prove that any\textbf{\ }$w\in
N_{U_{\varepsilon }}(u_{\varepsilon }^{\ast })$ can be written in the form
\begin{equation}
w=-\rho ^{\prime }(x)+\mu (x)\mbox{ a.e. }x\in (0,L)  \label{69}
\end{equation}%
with $\rho $ and $\mu $ previously defined.

Indeed, if $w$ is given by (\ref{69}) then
\begin{eqnarray*}
&&\int_{0}^{L}w(x)(u_{\varepsilon }^{\ast }-u_{\varepsilon
})dx=\int_{0}^{L}\mu (x)(u_{\varepsilon }^{\ast }-u_{\varepsilon
})dx-\int_{0}^{L}\rho ^{\prime }(x)(u_{\varepsilon }^{\ast }-u_{\varepsilon
})dx \\
&=&\int_{0}^{L}\mu (x)(u_{\varepsilon }^{\ast }-u_{\varepsilon })dx-\left.
\rho (u_{\varepsilon }^{\ast }-u_{\varepsilon })\right\vert
_{0}^{L}+\int_{0}^{L}\rho (x)(u_{\varepsilon x}^{\ast }-u_{\varepsilon
x})dx\geq 0
\end{eqnarray*}%
for any $u_{\varepsilon }\in U_{\varepsilon },$ by (\ref{70}) for the first
term and (\ref{44}) for the last term. This means that $w\in
N_{U_{\varepsilon }}(u_{\varepsilon }^{\ast }).$

For the inverse implication we take $w\in N_{U_{\varepsilon
}}(u_{\varepsilon }^{\ast })$ and we note first that $U_{\varepsilon }$ can
be written $U_{\varepsilon }=U_{1\varepsilon }+U_{2\varepsilon },$ where
\begin{equation*}
U_{1\varepsilon }=\{v\in W^{1,\infty }(0,L);\mbox{ }\left\vert
v_{x}\right\vert \leq u_{\infty }\mbox{ a.e. }x\in (0,L),\mbox{ }%
v(0)=u_{0}^{\varepsilon },\mbox{ }v(L)=u_{L}^{\varepsilon }\},
\end{equation*}%
\begin{equation*}
U_{2\varepsilon }=\{v\in L^{\infty }(0,L);\mbox{ }u_{m}(x)+\varepsilon \leq
v(x)\leq u_{M}(x)+2\varepsilon \mbox{ a.e. }x\in (0,L)\}.
\end{equation*}

We also remark that $U_{1\varepsilon }\cap $ int $U_{2\varepsilon }\neq
\varnothing ,$ and so
\begin{equation}
\partial (I_{1}+I_{2})=\partial I_{1}+\partial I_{2},  \label{69-2}
\end{equation}%
(see \cite{vb3}, p. 7), where $I_{i}$ are the indicator functions of $%
U_{i\varepsilon }$ ($i=1,2)$ and $\partial I_{i}$ denote their
subdifferentials.

Therefore, $w\in N_{U_{\varepsilon }}(u_{\varepsilon }^{\ast })=\partial
I_{U_{\varepsilon }}(u_{\varepsilon }^{\ast })$ is given by
\begin{equation}
w=\xi +\mu ,\mbox{ }\xi \in \partial I_{1}(u_{\varepsilon }^{\ast }),\mbox{ }%
\mu \in \partial I_{2}(u_{\varepsilon }^{\ast }).  \label{69-3}
\end{equation}

It is obvious that $\mu (x)\in \partial I_{2}(u_{\varepsilon }^{\ast
})=N_{U_{2\varepsilon }}(u_{\varepsilon }^{\ast })$ if and only if $\mu (x)$
satisfies (\ref{70}).

The form of $\xi =-\rho ^{\prime }$ (with $\rho $ given by (\ref{44}))
follows by some arguments based on general results given in \cite{vb3}, p.
11-15. For the reader's convenience we give a few details adapted to our
case.

Let $\xi \in \partial I_{1}(u_{\varepsilon }^{\ast })$ and $\gamma \in
W^{1,1}(0,L)$ such that $\xi =-\gamma ^{\prime }$ a.e. in $(0,L)$. This
means that
\begin{equation}
\xi =-\gamma ^{\prime }\in \partial I_{1}(u_{\varepsilon }^{\ast })
\label{42-0}
\end{equation}%
and $\gamma \in \partial I_{F}(u_{\varepsilon x}^{\ast }),$ where $I_{F}$ is
the indicator function of the set
\begin{equation*}
F=\{\zeta \in L^{\infty }(0,L);\mbox{ }\zeta =v^{\prime }\mbox{ a.e. in }%
(0,L),v\in U_{1\varepsilon }\}.
\end{equation*}%
Indeed, for all $v\in U_{1\varepsilon },$ we have%
\begin{equation*}
0\leq \int_{0}^{L}\xi (u_{\varepsilon }^{\ast }-v)dx=\int_{0}^{L}\gamma
(u_{\varepsilon x}^{\ast }-v_{x})dx=\int_{0}^{L}\gamma (u_{\varepsilon
x}^{\ast }-\zeta )dx,\mbox{ for all }\zeta \in F.
\end{equation*}%
The set $F$ can be decomposed as $F_{1}\cap F_{2}$ where
\begin{equation*}
F_{1}=\{\zeta \in L^{\infty }(0,L);\mbox{ }\zeta =v^{\prime },\mbox{ }v\in
W^{1,\infty }(0,L),\mbox{ }v(0)=u_{0}^{\varepsilon },\mbox{ }%
v(L)=u_{L}^{\varepsilon }\},
\end{equation*}%
\begin{equation*}
F_{2}=\{\zeta \in L^{\infty }(0,L);\mbox{ }\left\vert \zeta (x)\right\vert
\leq u_{\infty }\mbox{ a.e. }x\in (0,L)\}.
\end{equation*}%
Let us assume that there exists\textit{\ }$w_{0}\in W^{1,\infty }(0,L)$%
\textit{\ }such that\textit{\ }%
\begin{equation}
\left\Vert w_{0}\right\Vert _{L^{\infty }(0,1)}<u_{\infty },\mbox{ }%
w_{0}(0)=u_{0}^{\varepsilon },\mbox{ }w_{0}(L)=u_{L}^{\varepsilon }.
\label{42}
\end{equation}%
We note that by (\ref{42}), $w_{0}\in F_{1}\cap $ int $F_{2}$ and so
\begin{equation*}
\partial I_{F}=\partial I_{F_{1}}+\partial I_{F_{2}}
\end{equation*}%
(see again \cite{vb3}, p. 7). Therefore, $\gamma \in \partial
I_{F}(u_{\varepsilon x}^{\ast })$ can be written $\gamma =\gamma _{1}+\gamma
_{2}$ with $\gamma _{i}\in \partial I_{F_{i}}(u_{\varepsilon x}^{\ast }).$
The subdifferential $\partial I_{F}$ is an application from $L^{\infty
}(0,L) $ to $(L^{\infty }(0,L))^{\prime }$ and so $\gamma _{i}$ are seen as
elements belonging to $(L^{\infty }(0,L))^{\prime },$ being represented as
the sum of a continuous part $\gamma _{ia}\in L^{1}(0,L)$ and a singular
part $\gamma _{is}$ (see \cite{vb3}, p. 15). Then, $\gamma _{2a}\in \partial
I_{F_{2}}(u_{\varepsilon x}^{\ast })$ a.e. $x\in (0,L)$ and by Proposition
1.9 in \cite{vb3}, p. 11-13, we have
\begin{equation}
\gamma _{2a}(x)=\left\{
\begin{array}{l}
0\mbox{ \ \ \ \ \ \ \ \ \ \ \ \ \ a.e. in }\{x\in (0,L);\left\vert
u_{\varepsilon x}^{\ast }(x)\right\vert <u_{\infty }\} \\
\nu (x)u_{\varepsilon x}^{\ast }(x)\mbox{ \ a.e. in }\{x\in (0,L);\left\vert
u_{\varepsilon x}^{\ast }(x)\right\vert =u_{\infty }\},%
\end{array}%
\right.  \label{69-1}
\end{equation}%
where $\nu \in L^{1}(0,L),$ $\nu \geq 0$ a.e. in $(0,L).$

Now, $\gamma _{1a}\in \partial I_{F_{1}}(u_{\varepsilon x}^{\ast })$ and so%
\begin{equation*}
\int_{0}^{L}\gamma _{1a}(x)(u_{\varepsilon x}^{\ast
}-v_{x})(x)dx=-\int_{0}^{L}\gamma _{1a}^{\prime }(x)(u_{\varepsilon }^{\ast
}-v)(x)dx\geq 0,\mbox{ for any }v\in F_{1}.
\end{equation*}%
In particular setting $v:=u_{\varepsilon }^{\ast }+l\phi $ with $\phi \in
C_{0}^{\infty }(0,L)$ and $l>0,$ we get
\begin{equation*}
\int_{0}^{L}\gamma _{1a}^{\prime }(x)\phi (x)dx\geq 0\mbox{ for any }\phi
\in C_{0}^{\infty }(0,L).
\end{equation*}%
Setting $v:=u_{\varepsilon }^{\ast }-l\phi $ we get the inverse inequality
and so it follows that
\begin{equation*}
\int_{0}^{L}\gamma _{1a}^{\prime }(x)\phi (x)dx=0\mbox{ for any }\phi \in
C_{0}^{\infty }(0,L),
\end{equation*}%
which implies that $\gamma _{1a}^{\prime }(x)=0.$

In conclusion, $\gamma =\gamma _{1a}+\gamma _{2a},$ where $(\gamma
_{1a})^{\prime }=0$ a.e. in $(0,L)$ and $\gamma _{2}:=\rho $ satisfies (\ref%
{69-1}), so that by (\ref{42-0}) we have $\xi =-\gamma ^{\prime }=-\rho
^{\prime }.$ We note that
\begin{equation}
-\rho ^{\prime }\in \partial I_{1}(u_{\varepsilon }^{\ast })\mbox{ iff }\rho
\in \partial I_{F}(u_{\varepsilon x}^{\ast }).  \label{69-5}
\end{equation}%
By (\ref{69-3}) we get (\ref{69}) and thus relation (\ref{71}) follows from $%
\Phi (x)\in N_{U_{\varepsilon }}(u_{\varepsilon }^{\ast }).$

On the subset $\{x\in (0,L);$ $\left\vert u_{\varepsilon x}^{\ast
}(x)\right\vert =u_{\infty }$ a.e.$\}$ we have $\rho (x)=\nu
(x)u_{\varepsilon x}^{\ast }(x)$ and we observe that there are two cases for
$\nu >0:$%
\begin{equation*}
u_{\varepsilon x}^{\ast }(x)=\left\{
\begin{array}{c}
u_{\infty }\mbox{ \ if }u_{\varepsilon x}^{\ast }(x)=\frac{\rho (x)}{\nu (x)}%
>0 \\
-u_{\infty }\mbox{ if }u_{\varepsilon x}^{\ast }(x)=\frac{\rho (x)}{\nu (x)}%
<0.%
\end{array}%
\right.
\end{equation*}%
On the subset $\{x\in (0,L);$ $\left\vert u_{\varepsilon x}^{\ast
}(x)\right\vert <u_{\infty }$ a.e.$\}$ the function $\rho (x)=0.$ Therefore,
we deduce (\ref{72}), as claimed. \hfill $\square $

\section{Homogeneous Dirichlet-Neumann boundary conditions}

In this section we consider the final time minimization problem
\begin{equation}
\mbox{Minimize }\left( \int_{0}^{L}y^{u}(T,x)dx\right) \mbox{ for all }u\in
U,  \tag{$P_{1}$}
\end{equation}%
subject to
\begin{equation}
\frac{\partial y}{\partial t}-(uy_{x})_{x}=f\mbox{ in }Q,  \label{1p}
\end{equation}%
\begin{equation}
y(0,x)=y_{0}(x)\mbox{ in }(0,L),  \label{2p}
\end{equation}%
\begin{equation}
y(t,0)=0,\mbox{ }y_{x}(t,L)=0\mbox{ in }(0,T).  \label{3p}
\end{equation}%
The set $U$ and all hypotheses made for the functions occurring in $U$ are
the same as in Section 2.2.

The existence of the solution to the new state system is treated by the
variational technique, as in the case of the previous state system.

We introduce the space
\begin{multline*}
\widetilde{H_{u}^{1}}(0,L)=\{y\in L^{2}(0,L);\mbox{ }y\mbox{ locally
absolutely continuous in }[0,L]\backslash \{x_{0}\}, \\
\sqrt{u}y_{x}\in L^{2}(0,L),\mbox{ }y(0)=0\},
\end{multline*}%
equipped with the scalar product%
\begin{equation*}
\left( \theta ,\overline{\theta }\right) _{\widetilde{H_{u}^{1}}%
(0,L)}=\int_{0}^{L}\theta \overline{\theta }dx+\int_{0}^{L}u\theta _{x}%
\overline{\theta }_{x}dx,\mbox{ for any }\theta ,\overline{\theta }\in
\widetilde{H_{u}^{1}}(0,L).
\end{equation*}%
It is a Hilbert space with the norm given by (\ref{5}) and with the dual
denoted by $(\widetilde{H_{u}^{1}}(0,L))^{\prime }$. We have the continuous
and dense embeddings $\widetilde{H_{u}^{1}}(0,L)\hookrightarrow
L^{2}(0,L)\hookrightarrow (\widetilde{H_{u}^{1}}(0,L))^{\prime }$.

For the subscripts we denote $\widetilde{V_{u}}=\widetilde{H_{u}^{1}}(0,L),$
$\widetilde{V_{u}^{\prime }}=(\widetilde{H_{u}^{1}}(0,L))^{\prime }$ and $%
H=L^{2}(0,L).$

The solution to (\ref{1p})-(\ref{3p}) is defined according to Definition 2.3
by replacing $H_{u}^{1}(0,L)$ by $\widetilde{H_{u}^{1}}(0,L).$

Let $u\in U.$ We define $\widetilde{A}:\widetilde{H_{u}^{1}}(0,L)\rightarrow
(\widetilde{H_{u}^{1}}(0,L))^{\prime },$%
\begin{equation*}
\left\langle \widetilde{A}y,\psi \right\rangle _{\widetilde{V_{u}^{\prime }},%
\widetilde{V_{u}}}=\int_{0}^{L}uy_{x}\psi _{x}dx,\mbox{ for any }\psi \in
\widetilde{H_{u}^{1}}(0,L).
\end{equation*}%
Then, problem (\ref{1p})-(\ref{3p}) has a unique solution
\begin{equation}
y\in C([0,T];L^{2}(0,L))\cap L^{2}(0,T;\widetilde{H_{u}^{1}}(0,L))\cap
W^{1,2}([0,T];(\widetilde{H_{u}^{1}}(0,L))^{\prime })  \label{69-0}
\end{equation}%
obtained by the first part of Theorem 2.4.

\medskip

\noindent \textbf{Theorem 3.1. }\textit{Let }$y_{0}\in L^{2}(0,L)$\textit{\
and }$f\in L^{2}(Q),$ $y_{0}\geq 0$ \textit{a.e. on }$(0,L)$\textit{, }$%
f\geq 0$\textit{\ a.e. on} $Q.$\textit{\ Then }$(P_{1})$\textit{\ has at
least one solution.}

\medskip

The proof is led as in Theorem 2.5.

\subsection{Optimality conditions in the homogeneous Dirichlet-Neumann case}

\noindent \textbf{Proposition 3.2. }\textit{Let the assumptions of}
Proposition 2.9 \textit{hold. Assume in addition that }%
\begin{equation}
y_{0}\in H^{2}(0,L),\mbox{ }f\in C([0,T];H^{1}(0,L)),\mbox{ }y_{0x}>0\mathit{%
\ }\mbox{\textit{on}}\mathit{\ }[0,L],\mathit{\ }f_{x}>0\mathit{\ on\ }Q,
\label{75-0}
\end{equation}%
\begin{equation}
u_{M}\in W^{1,\infty }(0,L),\mbox{ }\left\Vert u_{M}^{\prime }\right\Vert
_{L^{\infty }(0,L)}\leq u_{\infty }.  \label{75}
\end{equation}%
\textit{Then, }$u^{\ast },$ \textit{a solution to }$(P_{1})$ \textit{has the
form }%
\begin{equation}
u^{\ast }(x)=\left\{
\begin{array}{l}
u_{0}+u_{\infty }x\mbox{ \ \ \ \ \ \ \ \ \ \ \ \ \textit{for} }x\in \lbrack
0,x_{1}) \\
u_{M}(x)\mbox{ \ \ \ \ \ \ \ \ \ \ \ \ \ \ \ \ \textit{for} }x\in \lbrack
x_{1},x_{2}) \\
-u_{\infty }(x-L)+u_{L}\mbox{ \ \textit{for} }x\in \lbrack x_{2},L],%
\end{array}%
\right.  \label{73-0}
\end{equation}%
\textit{where }$x_{1}$ \textit{and} $x_{2}$ \textit{are the solutions to }%
\begin{equation}
u_{M}(x_{1})=u_{\infty }x_{1}+u_{0},\mbox{ }u_{M}(x_{2})=-u_{\infty
}(x_{2}-L)+u_{L}.  \label{74-0}
\end{equation}%
\textit{The function }$u^{\ast }$\textit{\ is unique for fixed }$u_{M},$%
\textit{\ }$u_{m},$\textit{\ }$u_{0},$\textit{\ }$u_{L},$\textit{\ }$%
u_{\infty }.$

\medskip

\noindent \textbf{Proof.} Since, in particular, $y_{0}\in L^{2}(0,L)$ and $%
f\in C(\overline{Q})$ it follows by Theorem 3.1 that ($P_{1})$ has at least
a solution $u^{\ast }.$ We shall deduce the optimality condition for $%
(P_{1}) $ by passing to the limit in the approximating problem%
\begin{equation}
\mbox{Minimize }\left( \int_{0}^{L}y^{u}(T,x)dx\right) \mbox{ for all }u\in
U_{\varepsilon },  \tag{$P_{1\varepsilon }$}
\end{equation}%
subject to the state system (\ref{1p})-(\ref{3p}) with $U_{\varepsilon }$
given by (\ref{Ueps}).

All existence results and computations for the optimality condition in $%
(P_{1\varepsilon })$ are led similarly as for $(P_{\varepsilon }),$ but the
dual system and the optimality condition are slightly modified due to the
new cost functional and boundary conditions.

Let $(u_{\varepsilon }^{\ast },y_{\varepsilon }^{\ast })$ be a solution to $%
(P_{1\varepsilon }).$ The new dual system is%
\begin{equation}
\frac{\partial p_{\varepsilon }}{\partial t}+(u_{\varepsilon }^{\ast
}(x)p_{\varepsilon x})_{x}=0\mbox{ in }Q,  \label{51}
\end{equation}%
\begin{equation}
p_{\varepsilon }(T,x)=-1\mbox{ in }(0,L),  \label{52}
\end{equation}%
\begin{equation}
p_{\varepsilon }(t,0)=p_{\varepsilon x}(t,L)=0\mbox{ in }(0,T),  \label{53}
\end{equation}%
(here $\lambda _{1}=\lambda _{3}=0,$ $\lambda _{2}=1),$ and the optimality
condition
\begin{equation*}
\int_{0}^{T}\int_{0}^{L}(u_{\varepsilon }^{\ast }-u_{\varepsilon
})y_{\varepsilon x}^{\ast }p_{\varepsilon x}dxdt\leq 0,\mbox{ for all }%
u_{\varepsilon }\in U_{\varepsilon },
\end{equation*}%
is obtained by similar computations as for (\ref{40})$.$

\noindent Since $y_{0}\in H^{2}(0,L),$ $y_{\varepsilon }^{\ast }$ and $%
p_{\varepsilon x},$ solutions to nondegenerate equations, are more regular%
\begin{equation*}
y_{\varepsilon }^{\ast },\mbox{ }p_{\varepsilon }\in
C^{1}([0,T];L^{2}(0,L))\cap C([0,T];H^{2}(0,L)),
\end{equation*}%
and in particular $y_{\varepsilon x}^{\ast },$ $p_{\varepsilon x}\in
C([0,T]\times \lbrack 0,L]),$ $y_{\varepsilon xx}^{\ast },$ $p_{\varepsilon
xx}\in C([0,T];L^{2}(0,L)).$ Then, we can write
\begin{equation}
\int_{0}^{L}(u_{\varepsilon }^{\ast }-u_{\varepsilon })\left(
\int_{0}^{T}y_{\varepsilon x}^{\ast }p_{\varepsilon x}dt\right) dx\leq 0,%
\mbox{ for all }u_{\varepsilon }\in U_{\varepsilon }.  \label{50}
\end{equation}%
In addition, the convergence of $(P_{1\varepsilon })$ to $(P_{1})$ as $%
\varepsilon \rightarrow 0$ follows as in Theorem 2.8.

Moreover, according to Proposition 2.9 we have
\begin{equation}
\widetilde{\Phi }(x)=-\rho ^{\prime }(x)+\mu (x)\mbox{ a.e. }x\in (0,L),
\label{71-1}
\end{equation}%
\begin{equation}
u_{\varepsilon x}^{\ast }(x)\in u_{\infty }\mbox{sign(}\rho (x))\mbox{ a.e.%
\textit{\ }}x\in (0,L),  \label{72-1}
\end{equation}%
where
\begin{equation*}
\widetilde{\Phi }(x):=-\int_{0}^{T}y_{\varepsilon x}^{\ast }p_{\varepsilon
x}dt.
\end{equation*}%
\textbf{\ }

We continue by establishing the sign of $\widetilde{\Phi }(x).$ The solution
$y_{\varepsilon }^{\ast }$ is continuous and nonnegative on $[0,T]\times
\lbrack 0,L]$ and the minimum is attained on the boundary. Then by the
strong maximum principle it follows that $y_{\varepsilon }^{\ast }>0$ on $%
(0,T)\times (0,L)$ and $-y_{\varepsilon x}^{\ast }(t,0)<0.$ Let $%
z=u_{\varepsilon }^{\ast }y_{\varepsilon x}^{\ast }.$ Then, it satisfies
\begin{eqnarray*}
z_{t}-u_{\varepsilon }^{\ast }z_{xx}^{\ast } &=&u_{\varepsilon }^{\ast }f_{x}%
\mbox{ in }Q, \\
z(0,x) &=&u_{\varepsilon }^{\ast }y_{0x}\mbox{ in }(0,L), \\
z_{x}(t,0) &>&0,\mbox{ }z_{x}(t,L)=0\mbox{ in }(0,T).
\end{eqnarray*}%
This nondegenerate problem has a unique solution, continuous and positive on
$(0,T)\times (0,L)$, by the same strong maximum principle$,$ hence $%
y_{\varepsilon x}^{\ast }(t,x)>0,$ for any ($t,x)\in Q.$

From the dual problem we get $p_{\varepsilon }\leq 0$ on $[0,T]\times
\lbrack 0,L]$ and by a similar argument $-p_{\varepsilon x}(t,0)>0.$ If $%
\omega =u_{\varepsilon }^{\ast }p_{\varepsilon x}$ then $\omega $ satisfies
\begin{eqnarray*}
\omega _{t}+u_{\varepsilon }^{\ast }\omega _{xx}^{\ast } &=&0\mbox{ in }Q, \\
\omega (0,x) &=&0\mbox{ in }(0,L), \\
\omega _{x}(t,0) &<&0,\mbox{ }\omega _{x}(t,L)=0\mbox{ in }(0,T)
\end{eqnarray*}%
and it follows that $p_{\varepsilon x}(t,x)<0$ for any $(t,x)\in Q.$

In conclusion, $\widetilde{\Phi }(x)>0$ for all $x\in (0,L)$ and by (\ref%
{71-1}) we note that
\begin{equation*}
\rho ^{\prime }(x)=\mu (x)-\widetilde{\Phi }(x)
\end{equation*}%
preserves a negative sign only for $\mu (x)\leq 0,$ i.e., on the subset
\begin{equation*}
U_{-}=\{x\in (0,L);u_{\varepsilon }^{\ast }(x)<u_{M}(x)+2\varepsilon \}.
\end{equation*}%
It means that on this set $\rho (x)\neq 0$ except at most one point, and so
we get by (\ref{72-1}) that%
\begin{equation*}
\left\vert u_{\varepsilon x}^{\ast }(x)\right\vert =u_{\infty }\mbox{ a.e.
in }\{x\in (0,L);\mbox{ }u_{\varepsilon }^{\ast }(x)<u_{M}(x)+2\varepsilon
\},
\end{equation*}%
which is the eikonal equation. This solution cannot actually be uniquely
determined by the conditions $u_{\varepsilon }^{\ast }(0)=u_{0}^{\varepsilon
}$, $u_{\varepsilon }^{\ast }(L)=u_{L}^{\varepsilon },$ unless one observes,
following a similar argument as in \cite{vb-friedman}, that the function $%
u_{\varepsilon }^{\ast }$ is the maximal element of the set
\begin{equation*}
D=\{z\in W^{1,\infty }(0,L);\mbox{ }\left\vert z^{\prime }(x)\right\vert
\leq u_{\infty }\mbox{ a.e. }x\in (0,L),\mbox{ }z(x)\leq u_{\varepsilon
}^{\ast }\mbox{ }\forall x\in \partial U_{-}\}.
\end{equation*}%
Here, $\partial U_{-}$ is the boundary of $U_{-}.$ To this end, let $z\in D,$
and use (\ref{71-1}) to get on the right-hand that%
\begin{eqnarray*}
&&\int_{U_{-}}\rho ^{\prime }(x)(u_{\varepsilon }^{\ast }(x)-z(x))^{-}dx \\
&=&\int_{U_{-}}\mu (x)(u_{\varepsilon }^{\ast }(x)-z(x))^{-}dx-\int_{U_{-}}%
\widetilde{\Phi }(x)(u_{\varepsilon }^{\ast }(x)-z(x))^{-}dx\leq 0.
\end{eqnarray*}%
On the other hand
\begin{eqnarray*}
&&\int_{U_{-}}\rho ^{\prime }(x)(u_{\varepsilon }^{\ast }(x)-z(x))^{-}dx \\
&=&\left. \rho (x)(u_{\varepsilon }^{\ast }(x)-z(x))^{-}\right\vert
_{\partial U_{-}}+\int_{U_{-}}\rho (x)(u_{\varepsilon x}^{\ast
}(x)-z_{x}(x))dx\geq 0.
\end{eqnarray*}%
Since $\rho ^{\prime }<0$ on $U_{-}$, we necessarily obtain that $%
(u_{\varepsilon }^{\ast }(x)-z(x))^{-}=0,$ meaning that $z(x)\leq
u_{\varepsilon }^{\ast }(x)$ for any $z\in D.$

In conclusion, $u_{\varepsilon }^{\ast }(x)$ must have the slope equal to
either $u_{\infty }$ or $-u_{\infty }$ in the subset where $u_{\varepsilon
}^{\ast }(x)<u_{M}+2\varepsilon $, and being the maximal element of $D,$ it
is unique. In the subset where $u_{\varepsilon }^{\ast
}(x)=u_{M}(x)+2\varepsilon $ its slope is $u_{M}^{\prime }\in (-u_{\infty
},u_{\infty }).$ Therefore, for given data $u_{M},$ $u_{m},$ $%
u_{0}^{\varepsilon },$ $u_{L}^{\varepsilon },$ $u_{\infty },$ the possible
representation of $u_{\varepsilon }^{\ast }$ is%
\begin{equation}
u_{\varepsilon }^{\ast }(x)=\left\{
\begin{array}{l}
u_{0}^{\varepsilon }+u_{\infty }x\mbox{ \ \ \ \ \ \ \ \ \ \ \ \ \ for }x\in
\lbrack 0,x_{1}^{\varepsilon }) \\
u_{M}(x)+2\varepsilon \mbox{ \ \ \ \ \ \ \ \ \ \ \ for }x\in \lbrack
x_{1}^{\varepsilon },x_{2}^{\varepsilon }) \\
-u_{\infty }(x-L)+u_{L}^{\varepsilon }\mbox{ \ \ for }x\in \lbrack
x_{2}^{\varepsilon },L],%
\end{array}%
\right.  \label{73}
\end{equation}%
where $x_{1}^{\varepsilon }$\textit{\ }is the abscissa of the intersection
of the graphic of the function $x\rightarrow u_{M}(x)+2\varepsilon $ with
the line of slope $u_{\infty }$ passing through $u_{0}^{\varepsilon },$ and $%
x_{2}^{\varepsilon }$\ is the abscissa of the intersection of the graphic of
the function $x\rightarrow u_{M}(x)+2\varepsilon $ with the line of slope $%
u_{\infty }$ passing through $u_{L}^{\varepsilon },$ i.e.,%
\begin{equation}
u_{M}(x_{1}^{\varepsilon })+2\varepsilon =u_{\infty }x_{1}^{\varepsilon
}+u_{0}^{\varepsilon },\mbox{ }u_{M}(x_{2}^{\varepsilon })+2\varepsilon
=-u_{\infty }(x_{2}^{\varepsilon }-L)+u_{L}^{\varepsilon }.  \label{74}
\end{equation}

Passing to the limit in (\ref{73}) as $\varepsilon \rightarrow 0,$ on the
basis of the convergence of $(P_{1\varepsilon })$ to $(P_{1})$ we get (\ref%
{73-0}). \hfill $\square $

\medskip

As an example, the feature of a possible $u^{\ast }$ given by (\ref{73-0})
is seen in Fig. 3, by a solid line. It was computed for
\begin{eqnarray*}
L &=&1,\mbox{ }T=2,\mbox{ }x_{0}=0.7, \\
u_{M}(x) &=&20(x-x_{0})^{2},\mbox{ }u_{m}(x)=(x-x_{0})^{2}, \\
u_{0} &=&3,\mbox{ }u_{L}=1,\mbox{ }u_{\infty }=10.
\end{eqnarray*}%

The graphic of $u_{M}$ is drawn by a dashed line.
\begin{center}
        \begin{tabular}{c} 
           \includegraphics[width=3in]{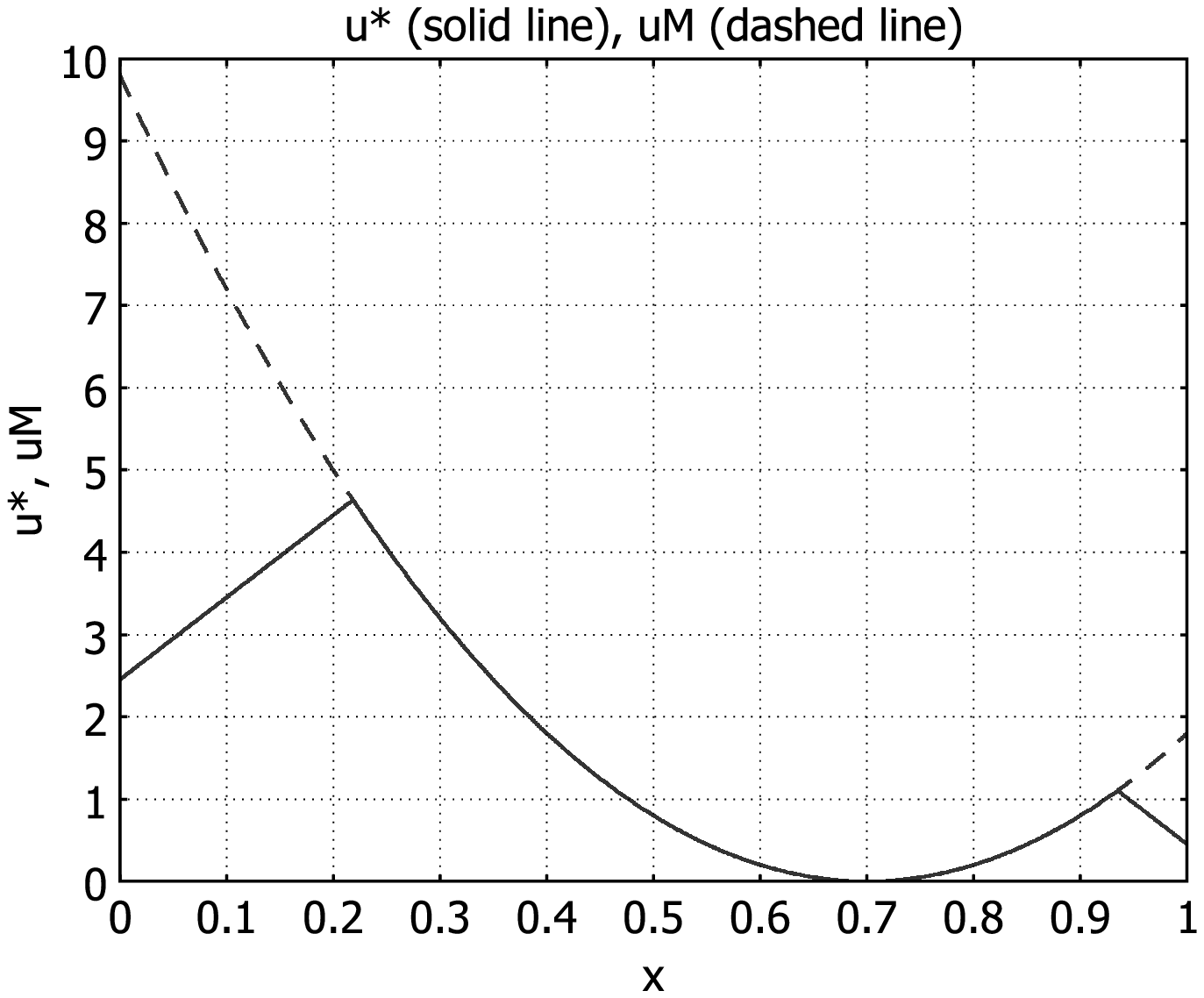} \\
        \end{tabular}

    Fig. 3. $u^{\ast }(x)$ given by (\ref{73-0}) and $u_{M}(x)$
    \end{center}

\section{Numerical results}

We note that the form of $u^{\ast }$ depends on the restrictions included in
the set $U,$ namely on the choice of $u_{M},$ $u_{m},$ $u_{\infty },$ $u_{0}$
and $u_{L},$ the greater importance being that of $u_{M}.$ We present some
numerical simulations to put into evidence the influence of $u_{M}$ in the
case when it is a polynomial function. The data are%
\begin{eqnarray*}
L &=&10,\mbox{ }T=0.5,\mbox{ }x_{0}=5,\mbox{ }y_{0}(x)=\sin (4x)+\left\vert
\sin (4x)\right\vert ,\mbox{ }f=0, \\
u_{M}(x) &=&2\left\vert x-x_{0}\right\vert ^{n},\mbox{ }u_{m}(x)=\left\vert
x-x_{0}\right\vert ^{n}, \\
u_{\infty } &=&10,\mbox{ }u_{0}=(u_{M}(0)+u_{m}(0))/2,\mbox{ }%
u_{L}=(u_{M}(L)+u_{m}(L))/2,\mbox{ }
\end{eqnarray*}%
and $n=3$ (Fig. 4.1, 4.2), $n=2$ (Fig. 5.1, 5.2), $n=1$ (Fig. 6.1, 6.2). The
figures represent the graphics of $u^{\ast }(x)$ computed by (\ref{73-0})
and the graphics of $y^{\ast }(T,x)$ representing the corresponding solution
to (\ref{1p})-(\ref{3p}) at $T=0.5.$ All computations are done with Comsol
Multiphysics v. 3.5a (FLN License 1025226).

\noindent The graphics of $y^{\ast }(T,x)$ in Figs. 4.2, 5.2, 6.2 and the
values $I=\int_{0}^{L}y^{\ast }(T,x)dx$ indicate $n=3$ as the best value
because $I$ is minimum in this case. The practice of the various physical
processes to which these results can be applied can indicate the choice of $%
u_{M}(x)$ and $u_{m}(x)$ in classes of functions other than polynomial, and
this may lead to better results.

\begin{center}
    \begin{tabular}{cc} 
           \includegraphics[width=2.6in]{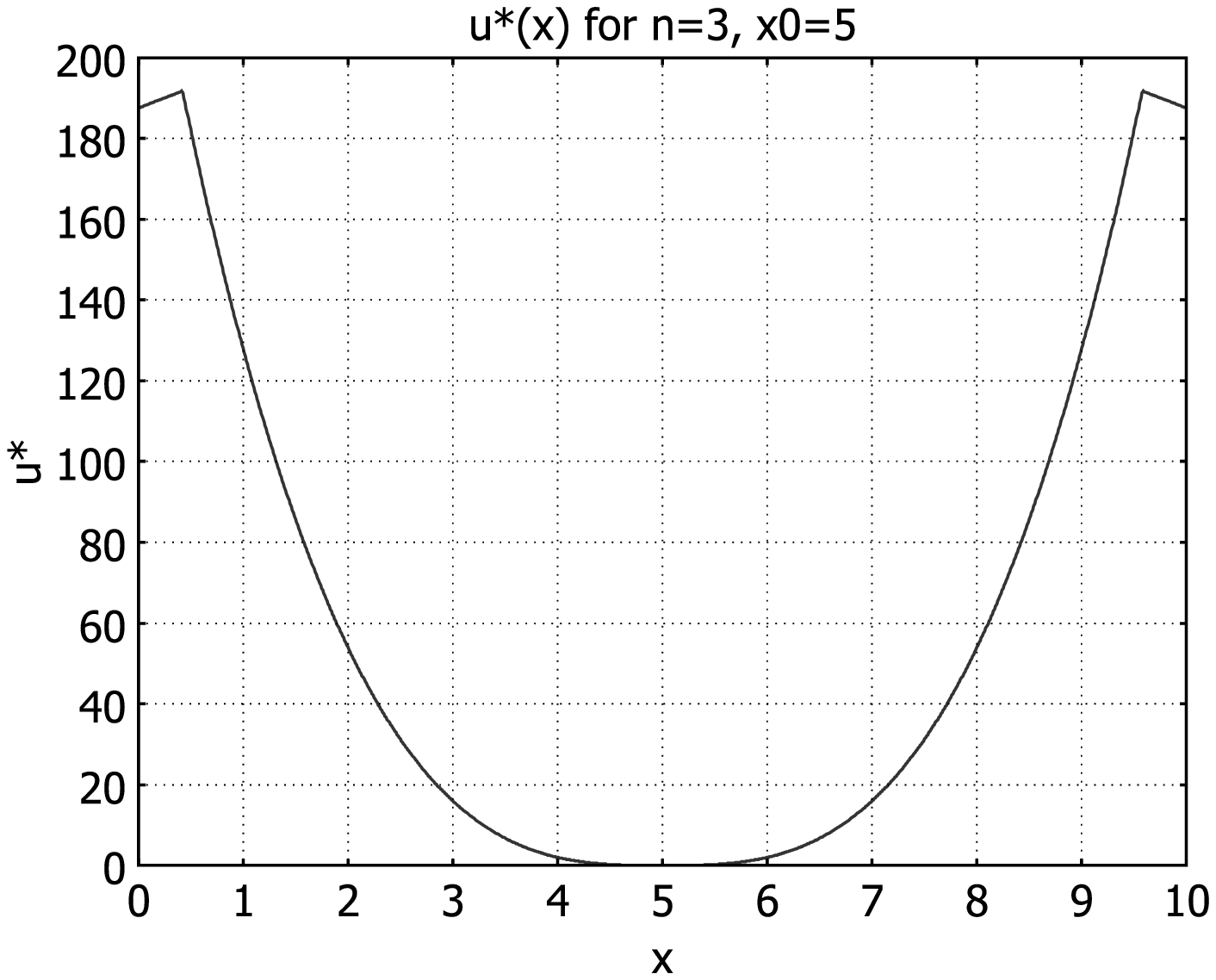} & \includegraphics[width=2.6in]{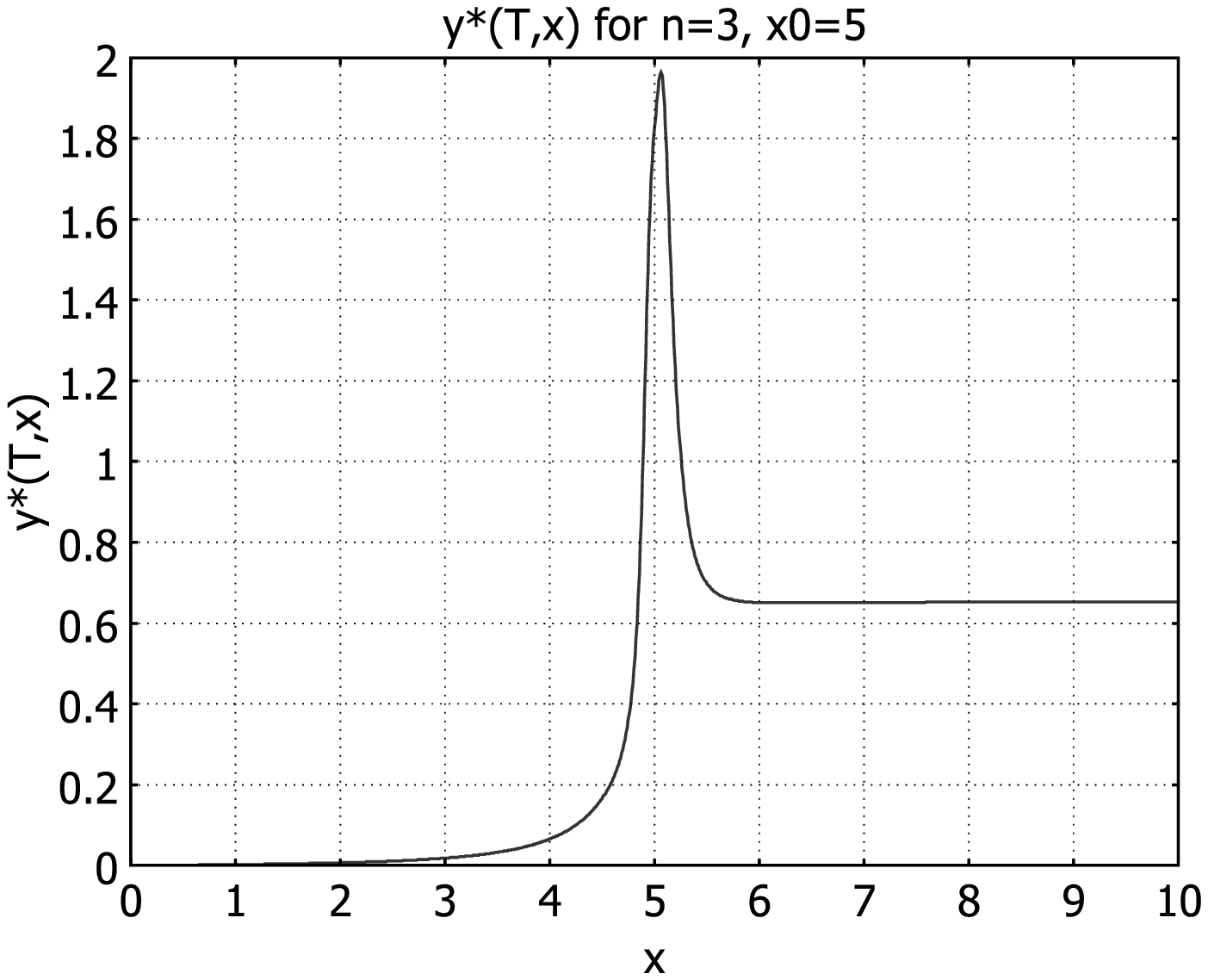} \\
             Fig. 4.1. $u^{\ast }(x)$ for $n=3$ & Fig. 4.2. $y^{\ast }(T,x)$ for $n=3$         \\ 
    \end{tabular}

\end{center}

\begin{center}
    \begin{tabular}{cc} 
           \includegraphics[width=2.6in]{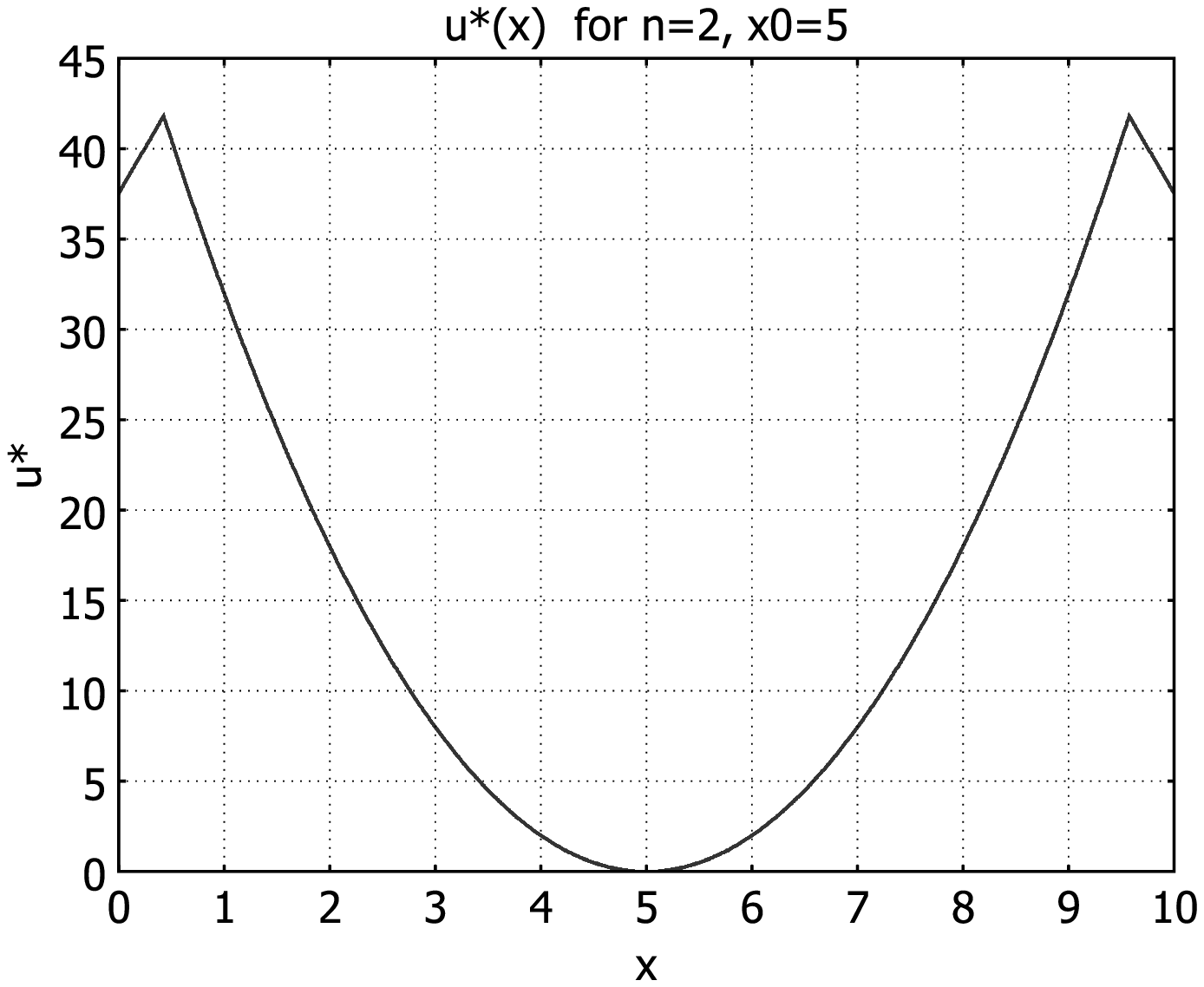} & \includegraphics[width=2.6in]{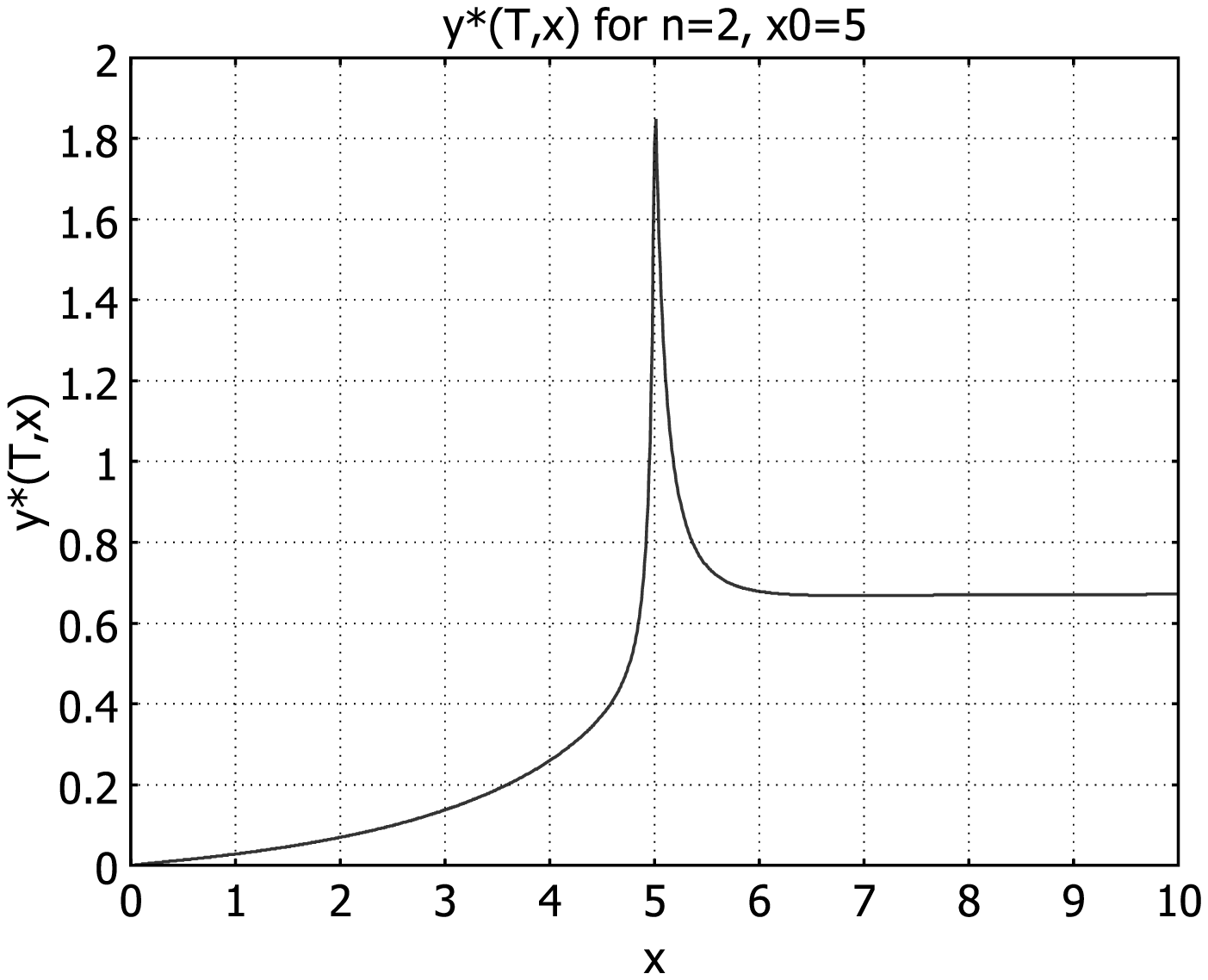} \\
           Fig. 5.1. $u^{\ast }(x)$ for $n=2$ & Fig. 5.2. $y^{\ast }(T,x)$ for $n=2$           \\ 
    \end{tabular}

\end{center}

\begin{center}
    \begin{tabular}{cc} 
           \includegraphics[width=2.6in]{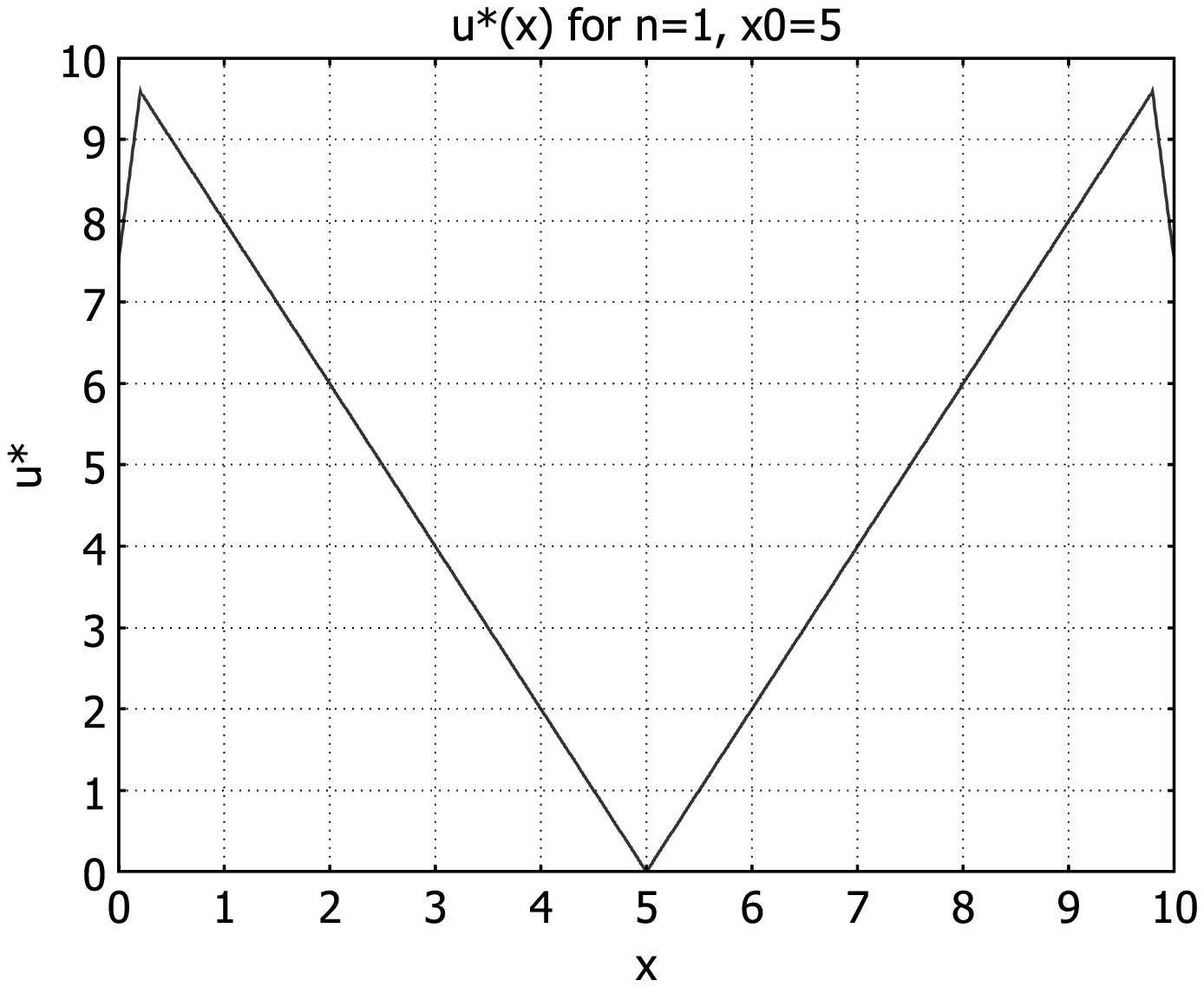} & \includegraphics[width=2.6in]{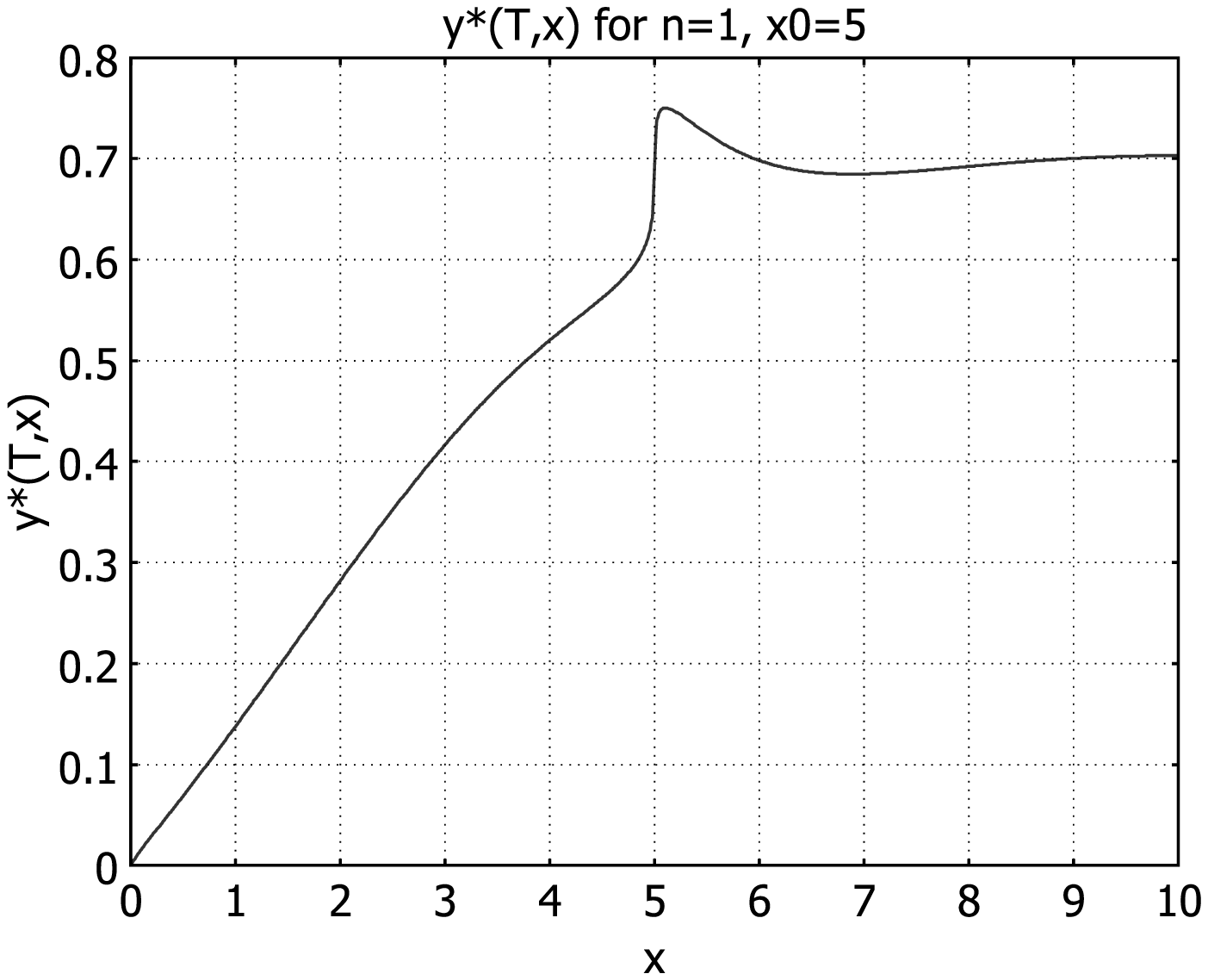} \\
             Fig. 6.1. $u^{\ast }(x)$ for $n=1$ & Fig. 6.2. $y^{\ast }(T,x)$ for $n=1$        \\ 
    \end{tabular}
\end{center}

\medskip

\noindent \textbf{Acknowledgments. }G.M. acknowledges the support\textbf{\ }%
of INDAM-GNAMPA, Italy, for May-June 2013 and of the grant CNCS--UEFISCDI,
project number PN-II-ID-PCE-2011-3-0027. G.F., R.M.M., S.R. acknowledge the
support\textbf{\ }of the GNAMPA project "Equazioni di evoluzione degeneri e
singolari: controllo e applicazioni", 2013.

\end{document}